\documentclass[CJK,12pt]{article}
\usepackage{amscd}
\usepackage{mathrsfs}
\usepackage{amsfonts}
\usepackage{graphicx}
\usepackage{bookmark}
\usepackage{amsmath,amscd,stmaryrd,amssymb,array, amsfonts,amsmath,amssymb}

\numberwithin{figure}{section}
 \numberwithin{equation}{section}

\newtheorem{theorem}{Theorem}[section]
\newtheorem{proposition}[theorem]{Proposition}
\newtheorem{definition}[theorem]{Definition}

\newtheorem{lemma}[theorem]{Lemma}
\newtheorem{remark}[theorem]{Remark}

\newcommand{\lb}{\label}

\newcommand{\ef}{\eqref}

\newcommand{\cA}{{\mathcal A}}
\newcommand{\cB}{{\mathcal B}}

\newcommand{\cC}{{\mathcal C}}

\newcommand{\cF}{{\mathcal F}}
\newcommand{\cH}{{\mathcal H}}
\newcommand{\cG}{{\mathcal G}}
\newcommand{\cO}{{\mathcal O}}

\newcommand{\cL}{{\mathcal L}}

\newcommand{\cZ}{{\mathcal Z}}
\newcommand{\cM}{{\mathcal M}}
\newcommand{\cN}{{\mathcal N}}
\newcommand{\cU}{{\mathcal U}}
\newcommand{\cV}{{\mathcal V}}
\newcommand{\cR}{{\mathcal R}}
\newcommand{\cS}{{\mathcal S}}
\newcommand{\cW}{{\mathcal W}}

\newcommand{\cK}{{\mathcal K}}

\newcommand{\sC}{{\mathscr C}}
\newcommand{\sE}{{\mathscr E}}
\newcommand{\sU}{{\mathscr U}}
\newcommand{\sM}{{\mathscr M}}
\newcommand{\sK}{{\mathscr K}}

\newcommand{\sX}{{\mathscr X}}
\newcommand{\sY}{{\mathscr Y}}

\newcommand{\sW}{{\mathscr W}}

\newcommand{\mB}{\mb{B}}

\def\be{\begin{equation}}
\def\ee{\end{equation}}
\def\ba{\begin{array}}
\def\ea{\end{array}}
\def\benu{\begin{enumerate}}
\def\eenu{\end{enumerate}}
\def\bt{\begin{theorem}}
\def\et{\end{theorem}}
\def\bp{\begin{proposition}}
\def\ep{\end{proposition}}
\def\bl{\begin{lemma}}
\def\el{\end{lemma}}
\def\br{\begin{remark}}
\def\er{\end{remark}}

\def\b{\beta}
\def\De{\Delta}
\def\de{\delta}
\def\pa{\partial}

\def\lam{\lambda}
\def\Lam{\Lambda}

\def\ve{\varepsilon}
\def\sig{\sigma}
\def\Sig{\Sigma}

\def\Om{\Omega}
\def\gam{\gamma}

\def\a{\alpha}

\def\.{\cdot}

\def\R{\mathbb{R}}

\def\A{\forall}
\def\ol{\overline}

\def\Cap{\bigcap}\def\Cup{\bigcup}

\def\ra{\rightarrow}

\def\~{\widetilde}
\def\8{\infty}
\def\X{\times}
\def\({\left(}
\def\){\right)}
\def\[{\left[}
\def\]{\right]}

\def\mb{\mbox}
\def\emp{\emptyset}
\def\sm{\setminus}
\def\ss{\Subset}
\def\<{\langle}
\def\>{\rangle}

\def\Hs{\hspace{1cm}}\def\hs{\hspace{0.5cm}}
\def\Vs{\vskip8pt}\def\vs{\vskip4pt}

\def\({\left(}\def\){\right)}

\begin{document}

\begin{center}
{\bf\Large Equilibrium Index of Invariant Sets and Global\\Static Bifurcation for  Nonlinear\\[1ex] Evolution Equations}
\end{center}
\vs\centerline{Desheng  Li \footnote{Supported by the grant of NSFC (10771159,
11071185)}}  %\vskip20pt
\begin{center}
{\footnotesize {Department of Mathematics,  Tianjin University
%\\     and\\     Center  of Applied Mathematics,  Tianjin University
\\     Tianjin 300072,  China\\

{\em E-mail}:  lidsmath@tju.edu.cn}}
\end{center}

\centerline{Zhi-Qiang Wang\footnote{Corresponding author, supported by the grant of NSFC (11271201)}}

  %\vskip20pt
\begin{center}
{\footnotesize
{%Center  of Applied Mathematics,  Tianjin University\\     Tianjin 300072,  China\\and\\
     Department of Mathematics and Statistics\\
Utah State University Logan, UT 84322\\
 {\em E-mail}: zhi-qiang.wang@usu.edu}}
\end{center}
\Vs

{\bf Abstract.}  We introduce the notion of equilibrium index for statically isolated invariant sets of the system $u_t+A u=f_\lam(u)$ on Banach space $X$ (where $A$ is a sectorial operator with compact resolvent) and  present  a reduction theorem and an index formula for bifurcating invariant sets near equilibrium points.
Then we prove a new global static bifurcation theorem where the crossing number $\mathfrak{m}$ may be even. 
In particular, in case   $\mathfrak{m}=2$, we show that  the system  undergoes either an attractor/repeller bifurcation, or a global static bifurcation. An illustrating example is also given by considering the bifurcations of the  periodic boundary value problem of second-order differential equations.

%\newpage

\section{Introduction}

In our previous work \cite{LW} we studied the dynamic bifurcation of the  equation
\be\label{e:1.1}
u_t+A u=f_\lam(u)
\ee on a Banach space $X$, where $\lam\in\R$ is the bifurcation parameter,   $A$ is a sectorial operator on $X$ with compact resolvent, and  $f_\lam(u)$ is a continuous  mapping from $X^\alpha\X \R$ to $X$ for some $0\leq\alpha<1$.
Suppose $f_\lam(0)\equiv0$ for $\lam\in\R,$ and hence $u=0$ is always a trivial equilibrium solution of  (\ref{e:1.1}). It was shown that if the crossing number 
%(the number of eigenvalues   of the linearized equation crossing the imaginary axis)
 at $\lam=\lam_0$ is nonzero  (the ``odd-multiplicity condition'' is not needed), then the system bifurcates from the equilibrium solution $0$ an isolated compact invariant set $K_\lam$ with nontrivial Conley index. Moreover, such a bifurcation enjoys some global features as in the classical Rabinowitz's global static bifurcation theorem.

Now a  natural problem arises: Does the bifurcated invariant set $K_\lam$ contain equilibrium solutions\,? For gradient-like systems, this question  seems to be somewhat  trivial, as any nonempty  compact invariant set of such a system necessarily contains an equilibrium.  Another particular but important case is the one of the attractor bifurcation, for which  Ma and Wang established  some index formulas via indices of isolated singular points (see e.g. \cite{MW1}). When the crossing number is odd, by using Ma and Wang's index formulas one can  assure the existence of equilibrium solutions in $K_\lam$.
%hence recover the classical Krasnosel'skii's Bifurcation Theorem.
Here we are interested in the general case. Our strategy is as follows.

First, we introduce the notion of equilibrium index for statically isolated invariant sets  for the non-parameterized equation \be\label{e1.2}u_t+A u=f(u)\ee on $X$, where $A$ is the same as in \eqref{e:1.1}, and  $f$ is a locally Lipschitz continuous  mapping from $X^\alpha$ to $X$ ($0\leq\alpha<1$).
%Here $X^\a=D((aI+A)^\alpha)$, and $a$ is   a real number such that $\mb{Re}\,\sig(aI+A)>0$.
Denote $\Phi$ the local semiflow generated by \eqref{e1.2}, and let $K$ be a statically isolated invariant set.
We define the equilibrium index $\mb{Ind}\,(\Phi,K)$ of an isolated invariant set $K$ to be the Leray-Schauder degree $\mb{deg}\,(I-F,\Om,0)$, where
$$F=(aI+A)^{-1}(aI+f)$$ which maps $X^\a$ into itself, $a$ is a real number such that $\mb{Re}\,\sig(aI+A)>0$, and $\Om$ is a statically isolating neighborhood of $K$ in $X^\a$. It can be shown that the index  $\mb{Ind}\,(\Phi,K)$ is independent of the choice of the number $a$ and  the neighborhood  $\Om$, and hence is well defined. One of the advantages of introducing this notion is that, it allows us to obtain information on equilibrium solutions by directly performing some mathematical  analysis on the evolution equation \eqref{e1.2} without returning  back to the corresponding stationary equation and putting it into a proper form so that one can apply the  Leray-Schauder degree. % For a  dynamically isolated invariant set $K$ of a finite dimensional system $u_t=f(u)$, it can be shown that $\mb{Ind}\,(\Phi,K)$ equals precisely the Euler number $\chi(h(\Phi,K))$ of the Conley index $h(\Phi,K)$.

Secondly, we prove a reduction theorem for equilibrium index defined above %  near equilibrium points
by using the geometric theory of evolution equations, which allows us to compute the index $\mb{Ind}\,(\Phi,K)$ of an invariant set $K$ near an equilibrium point $e$ by restring the system $\Phi$ on the local center manifold of $e$. More specifically, let $e$ be an equilibrium of $\Phi$. We show that there is a neighborhood $U$ of $e$ such that for any statically isolated invariant set $K$ of $\Phi$ in $U$, it holds that
$$
\mb{Ind}\,(\Phi,K)=(-1)^{\mathfrak{m}_1}\mb{Ind}\,(\Phi^c,K),
$$
where $\mathfrak{m}_1$ is the dimension of the local unstable manifold of $e$, and $\Phi^c$ is the restriction of $\Phi$ on the local center manifold.
 Since invariant manifolds are actually of a pure dynamical nature, the above result indicates that there are inherent connections between static and dynamic objects, which fact was recognized as early as in the work of Chow and Hale \cite{Chow}.

Thirdly, based on the reduction theorem and a result on the relation between topological degree and Conley index for finite dimensional systems given in Rybakowski \cite{Ryba}, we establish
an index formula for the bifurcating invariant set  of \eqref{e:1.1}. Denote $\Phi_\lam$ the local semiflow generated by \eqref{e:1.1}, and let $\Phi^c_\lam$ be the restriction of $\Phi_\lam$ on the local center manifold of $0$. Suppose $\lam=\lam_0$ is  a dynamic bifurcation value, and let $K_\lam$ be the isolated  invariant set  bifurcated from the trivial equilibrium $0$. We prove that
\be\label{e1.4}\ba{ll}
\mb{Ind}\,(\Phi_\lam,K_\lam)&=\left\{\ba{ll}\chi\(h(\Phi_{\lam_0},\{0\})\)-(-1)^{\mathfrak{m}_1},\hs&\lam<\lam_0;\\[1ex]
\chi\(h(\Phi_{\lam_0},\{0\})\)-(-1)^{\mathfrak{m}_1+\mathfrak{m}_2},\hs&\lam>\lam_0.\ea\right.
\\[4ex]
&=(-1)^{\mathfrak{m}_1}\left\{\ba{ll}\chi\(h(\Phi_{\lam_0}^{c},\{0\})\)-1,\hs&\lam<\lam_0;\\[1ex]
\chi\(h(\Phi_{\lam_0}^{c},\{0\})\)-(-1)^{\mathfrak{m}_2},\hs&\lam>\lam_0,\ea\right.
\ea
\ee
where $\mathfrak{m}_1$ and $\mathfrak{m}_2$ denote, respectively, the dimensions of the local  unstable manifold and the local center manifold of $0$ at $\lam=\lam_0$ ($\mathfrak{m}_2$ is precisely the crossing number), and $\chi(h(\Phi,K))$ denotes the Euler number  of the Conley index $h(\Phi,K)$.
Due to the degeneracy of the equilibrium $0$ at $\lam=\lam_0$, the computation of the index $\mb{Ind}\,(\Phi_{\lam_0},\{0\})$ via topological degree seems to be of little hope. However, in some cases  computing the  Conley indices  $h(\Phi_{\lam_0},\{0\})$ and $h(\Phi^c_{\lam_0},\{0\})$ may be of practical sense. For instance, in the particular case of attractor bifurcation, $0$ is an attractor of the $\Phi^c_{\lam_0}$. Hence one easily deduces that
$$h(\Phi^c_{\lam_0},\{0\})=\Sigma^0,
$$
where $\Sigma^0$ denotes the homotopy type of the pointed $0$-dimensional sphere (the space consisting of precisely two distinct points with one being the base point). Consequently $\chi\(h(\Phi_{\lam_0}^{c},\{0\})\)=1$. Thus for $\lam>\lam_0$, by \eqref{e1.4} we see that
\be\label{e1.5}\ba{lcl}
\mb{Ind}\,(\Phi_\lam,K_\lam)&=(-1)^{\mathfrak{m}_1}\left\{\ba{ll}0,\hs &\mb{$\mathfrak{m}_2$\,=\,even};\\[1ex]
2,\Hs &\mb{$\mathfrak{m}_2$\,=\,odd},\ea\right.
\ea
\ee
which, in our situation,  recovers an   index formula in the theory of attractor bifurcation given be Ma and Wang; see e.g. \cite[Theorem 6.1]{MW0}.
Another important example is the one where $\mathfrak{m}_2=2$ and the equilibrium $0$ fails to be an attractor of $\Phi_{\lam_0}^c$; see Section 8 for details.

Fourthly, we are interested in the global static bifurcation of \eqref{e:1.1}. A well-known result in this line is the famous Rabinowitz's Global Bifurcation Theorem. However, it requires the crossing number to be odd (``crossing odd-multiplicity'' condition).  If one drops  this  condition then situations become very complicated. To the authors' knowledge, even if for gradient systems the global static bifurcation still remains an open problem.
 %although by a classical bifurcation theorem on potential operator equations due to  Krasnosel'skii (see e.g. \cite[Chap.\,II, Sect.\,7]{Kie}) we know that static bifurcation occurs as long as the crossing number $\mathfrak{m}_2>0$,
 % in case  there are a finite number of eigenvalues of the linearized equation crossing the imaginary axis
% Whereas  the global bifurcation remains an open problem.
To deal with this problem without assuming the ``crossing odd-multiplicity'' condition, Ma and Wang \cite{MW4} proved some new local and global static  bifurcation theorems by using higher-order nondegenerate singularities of nonlinearities.  In this present work, motivated by the index formula \eqref{e1.4} we give a new global static bifurcation theorem. Roughly speaking, we show that if $\lam_0$ is a bifurcation value and
\be\label{e1.6}
\chi\(h(\Phi_{\lam_0}^{c},\{0\})\)\ne 1\mb{ or } (-1)^{\mathfrak{m}_2},
\ee
then the equation \eqref{e:1.1} bifurcates from the trivial stationary solution $(0,\lam_0)$ a connected branch $\Gamma$ of stationary solutions. (We call $(e,\lam)$ a stationary  solution of \eqref{e:1.1} if $e$ is an equilibrium point of $\Phi_\lam$.) $\Gamma$ enjoys some global features as in the classical Rabinowitz's Global Bifurcation Theorem. It is worth noticing that if the crossing number  $\mathfrak{m}_2$ is odd, then the condition \eqref{e1.6} is automatically satisfied.

 Finally, we pay some  special attention  to the case where $\mathfrak{m}_2=2$, i.e.,  there are exactly two eigenvalues (with multiplicity) crossing the imaginary axis. Let $(0,\lam_0)$ be a bifurcation point. Suppose $S_0=\{0\}$ is dynamically isolated with respect to the flow $\Phi_{\lam_0}$.  We prove that the system  either undergoes an attractor/repeller bifurcation (a generalized Hopf bifurcation), or bifurcates from $(0,\lam_0)$  a connected global  bifurcation branch of stationary solutions. As an illustrating example, we consider the bifurcation of the  periodic boundary value problem of the second-order differential equation \be\label{e1.7}-u''=\lam u+a(x)u^2+h(x,u),\Hs x\in\R,\ee where $h(x,u)=0(|u|^3)$ at $u=0$.
Let $A$ be the differential operator $-\frac{d^2}{dx^2}$ associated with  periodic boundary condition. Then
$$
\sig(A)=\sig_p(A)=\{\lam_n\}_{n=0}^\8.
$$
The first eigenvalue $\lam_0$ is simple with a corresponding constant eigenfunction, and all the others are of multiplicity $2$. For each $\lam_k$ ($k\geq 1$),
we show under appropriate conditions that either there is a two-sided neighborhood $I$ of $\lam_k$ such that for each $\lam\in I\sm\{0\}$, the problem has at least two distinct nontrivial solutions, or it bifurcates from $(0,\lam_k)$ a global bifurcation branch enjoying the properties in the Rabinowitz's Global Bifurcation Theorem.

This paper is organised as follows. In Section 2 we make some preliminaries,
and in Section 3 we introduce the notion of equilibrium index and discuss its basic properties. In Section 4 we establish a reduction theorem for equilibrium index. In Section 5 we give an index formula for the bifurcating invariant sets of \eqref{e:1.1}. Section 6 is devoted to the global static bifurcation of \eqref{e:1.1}, in which we prove a global static bifurcation theorem without the ``odd-multiplicity condition''. Section 7 consists of some discussions on the special case where the crossing number $\mathfrak{m}_2=2$. In Section 8 we give an illustrating  example by considering the periodic boundary value problem of \eqref{e1.7}.

\section{Preliminaries}

This section is concerned with some preliminaries.
\subsection{Basic topological notions and facts}
Let $X$ be  a  metric space with metric $d(\cdot,\cdot)$.
For convenience  we will always identify a singleton $\{x\}$ with the
 point $x$ for any $x\in X$.

Let $A$  be a nonempty subset of $X$. 
 The closure, interior and boundary of $A$ are denoted, respectively, by $\ol A$, int$\,A$ and $\pa A$. 
  A subset $U$ of $X$ is called a {\em neighborhood} of $A$ if $\ol A\subset \mbox{int}\,U$. The {$\ve$-neighborhood} of $A$ in $X$, noted by $\mB_X(A,\ve)$, is defined to be the set $\{y\in
X:\,d(y,A)<\ve\}.$

Let $A$ and $B$ be two nonempty subsets of $X$. 
If $A\subset B$ then  we will use the notations int$_B\,A$ and $\pa_B A$ to denote the interior and  boundary of $A$ in $B$, respectively.
The verification of the following basic facts is straightforward. 
\bl\label{bl} $\pa_B(A\cap B)\subset\pa A$. If $A\subset \mb{\em int}B$ then we also have $\pa_BA=\pa A$.
\el

The {\em distance} $d(A,B)$ between $A$ and $B$ is defined as
 $$
 d(A,B)=\inf\{d(x,y):\,\,\,x\in A,\,\,y\in B\}.
 $$
 The {\em Hausdorff semi-distance} and {\em  distance} of $A$ and $B$ are defined
as
$$d_{\mbox{\tiny H}}(A,B)=\sup_{x \in A}d(x,B),\hs
\delta_{\mbox{\tiny H}}(A,B) = \max\left\{d_{\mbox{\tiny H}}(A,B),
d_{\mbox{\tiny H}}(B,A)\right\}
$$
 respectively.
 We also assign
$d_{\tiny\mb{H}}(\emp,B)=0$.

Let $A_\lam$ ($\lam\in\Lam$) be a family of nonempty subsets of $X$, where $\Lam$ is a metric space. We say that $A_\lam$ is {\em upper semicontinuous} in $\lam$ at $\lam_0\in\Lam$, this means  $$d_{\tiny\mb{H}}(A_\lam,A_{\lam_0})\ra 0,\hs \mb{as }\,\lam\ra\lam_0.$$

The following two lemmas will play  important roles in our discussion.

\bl\label{l:2.2}\cite{Rab} Let $X$ be a compact metric space, and let $A$, $B$ be two disjoint closed subsets of $X$. Then either there is a subcontinuum $C$ of $X$ such that
$$\ba{ll}
A\cap C\ne \emp\ne B\cap C,\ea
$$
or $X=X_A\cup X_B$, where $X_A$ and $X_B$ are disjoint compact subsets of $X$ containing $A$ and $B$, respectively.

\el

\bl\label{l:2.3} \cite{CV} (pp. 41) Let $X$ be a compact metric space. Denote by $\sK(X)$ the family of compact subsets of $X$ which is equipped with the Hausdorff metric $\de_{\mbox{\tiny H}}(\.,\.)$. Then $\sK(X)$ is a compact metric space.
\el

%\subsection{Homotopy of topological spaces}
%\subsubsection{Homotopy of topological  spaces}
%Let $X$ and $Y$ be two topological spaces.

%We say that the mappings  $f,g:X\ra Y$  are {\em homotopic}, notated as $f\simeq g$, if there is a mapping $F:X\X[0,1]\ra Y$ such that $$F(\.,0)=f,\hs F(\.,1)=g.$$

%If there exist a  mapping $f:X\ra Y$ and a mapping $g:Y\ra X$ such that
%$$
%g\circ f\simeq id_X,\hs f\circ g\simeq id_Y,
%$$
%then we say that $X$ and $Y$ are {\em homotopy equivalent}, notated as $X\simeq Y$.

%For topological space $X$, if we identify a point $x_0\in X$ as the {\em base point}, then the  pair $(X,x_0)$ is called a {\bf pointed space}.
%One can also define homotopy equivalence as above for pointed spaces. The only difference is that in this case  one uses pointed mappings instead of the usual ones.  We omit the details.

%Let  $I =[0,1]$, and let $A$ be a closed subset of $X$. $A$ is said to be a {\bf strong deformation retract} of $X$, if there exists a mapping $H:X\X I\ra X$ such that $$H(\.,0)=id_X,\hs H(X,1)\subset A,$$ and
%$$
%H(x,\sig)=x,\Hs\mb{for all }x\in A.
%$$
%In this subsection we  present some fundamental criteria  on homotopy equivalence of spaces and collect some basic operations between pointed spaces.

\subsection{Some fundamental  dynamical notions}
In this subsection we collect some fundamental dynamical notions for the reader's convenience.

Let $X$ be a metric space.
\begin{definition}
A {\em local semiflow} $\Phi=\Phi(t,x)$ on $X$ is a continuous mapping from an
open set $\mathcal{D}(\Phi)\subset \mathbb{R}^+\times X$ to $X$ that
enjoys the following properties: \benu
\item[{\rm(1)}] for each $x\in X$, there exists $0<T_x\leq\infty$ such that
$$(t,x)\in\mathcal{D}(\Phi)\Longleftrightarrow t\in [0,T_x);$$

\item[{\rm(2)}] $\Phi(0,\cdot)=id_X$, and
$\Phi(t+s,x)=\Phi(t,\Phi(s,x))$
for all $x\in X$ and $t,s\in\R^+$ with $t+s\leq T_x$\,.
 \eenu
%\end{defn}
\end{definition}

Assume that there has been given a local semiflow $\Phi$ on $X$. As usual  we will write $\Phi(t,x)=\Phi(t)x$.
% (\ref{e2.1}) reads: $$ \Phi(s+t)x=\Phi(t)\Phi(s)x. $$

Let $J$ be an interval.  A {\it trajectory} (or {\em solution}) of $\Phi$ on $J$
is a continuous mapping $\gam:J\ra X$ such that
$$
  \gam(t)=\Phi(t-s)\gam(s),\hs \forall t,s\in J,\, t\geq s.
  $$
A trajectory $\gam$ on $J=\R$ is called  a {\it full trajectory}.
% A constant full trajectory $\gam(t)\equiv e$ is called an {\em equilibrium}.

%The {\em orbit} of a trajectory $\gam$ on $I$ is the set $$ \mb{orb}(\gam)=\{\gam(t):\,\,t\in I\}.$$ The orbit of a full trajectory is simply called a {\em full orbit}.

% The {\em $\omega$-limit set} $\omega(\gam)$ % and $\omega^*$-limit set  of a full trajectory $\gam$ is defined as $$\ba{ll}\omega(\gamma)=\{y\in X:\,\,\,\,\mb{there exists }  t_n\ra \8 \mb{ such that }\gamma(t_n)\ra y\}.\ea$$
%$$\ba{ll}\omega^*(\gamma)=\{y\in X:\,\,\,\,\mb{there exists }  t_n\ra -\8 \mb{ such that }\gamma(t_n)\ra y\}.\ea$$

%Given $U\subset X$,  denote $K_\8(\Phi,U)$ the {\em union of all bounded full orbits} in $U$. In the case where $U=X$, we will simply write $$K_\8(\Phi,X)=K_\8(\Phi).$$
% It is trivial to see that $K_\8(\Phi,U)$  is invariant.

%\subsubsection{Asymptotic compactness of semiflows}

A set $S\subset X$ is said to be   {\em positively invariant} (resp. {\em invariant}), if $\Phi(t)S\subset S$ (resp. $\Phi(t)S=S$) for all $t\geq0$.

A compact invariant set $\cA\subset X$ is said to be an {\em attractor} of  $\Phi$, if it attracts a neighborhood $U$ of itself, namely,  %$\Phi(t)x$ exists on $[0,\8)$ for all $x\in U$; moreover,
$
\lim_{t\ra\8}d_H\(\Phi(t)U,\cA\)=0.
$

The {\em attraction basin} of an attractor $\cA$, denoted by $\sU(\cA)$, is defined
as
$$\sU(\cA)=\{x:\,\,\lim_{t\ra\8}d(\Phi(t)x,\,\cA)=0\}.$$
As in \cite[Proposition 3.4]{Lijmaa}, one can easily verify  that $\sU(\cA)$ is open.

Let   $S$ be  a compact invariant set. Then  the restriction $\Phi_S$ is a semiflow on $S$.
A compact set $\cA\subset S$ is said to be an {\em attractor of $\Phi$ in $S$},
this  means that $\cA$ is an attractor of
$\Phi_S$ in $S$.
%\br
%In the case when $S$ is an attractor, any attractor  of $\Phi$ in $S$ is also an attractor in $X$; see, e.g., \cite{Kap}.
%\er

Given an attractor $\cA$
 of $\Phi$ in $S$, define
\be\label{e:3.0a}\ba{ll}\cA^*=\{x\in S:\,\, \omega(x)\Cap \cA=\emptyset\}.\ea\ee
 $\cA^*$ is called the {\em repeller }of
$\Phi$ in $S$  dual to $\cA$, and $(\cA,\cA^*)$ is called an {\em attractor-repeller pair} in $S$.

%Let $N\subset X$. %We say that $\Phi$ {\em does not explode} in $N$, if $$\Phi([0,T_x))x\subset N\Longrightarrow T_x=\infty.$$
%\begin{definition}\label{defn2.3}\cite{Ryba}\,
A subset $N$ of $X$ is said to be {\em admissible} \cite{Ryba}, if for any sequences $x_n\in
N$ and $t_n\rightarrow+\infty$ with $\Phi([0,t_n])x_n\subset N$, the sequence $\Phi(t_n)x_n$ has a convergent subsequence.
$N$ is said to be {\em strongly admissible} \cite{Ryba} if in addition, $\Phi$ does not explode in $N$, namely,
$$\Phi([0,T_x))x\subset N\Longrightarrow T_x=+\infty.$$
%\end{definition}

\begin{definition}\label{defn2.4}
$\Phi$ is said to be asymptotically compact  on $X$, if each bounded
set $B\subset X$ is strongly admissible.
\end{definition}

\br\label{r2.1}If $\Phi$ is asymptotically compact, then  one easily verifies that each bounded invariant set of $\Phi$ is necessarily precompact.
\er

\subsection{Conley index}
%In this subsection let us recall  briefly the definition of  Conley index.
From now on  we assume that $X$ is a {\em complete metric space}.
Although we don't require $X$ to be complete in the definition of a local semiflow, completeness of the phase  space always plays a crucial role in establishing a dynamical systems theory.

Let  $\Phi$ be  a local semiflow on $X$. Since we are working in an infinite dimensional space, in the remaining part of this section, we also assume that  %$\Phi$ satisfies the following {\em asymptotic compactness} assumption:
\begin{enumerate}
\item[({\bf AC})] $\Phi$ is asymptotically compact.
\end{enumerate}

  %Therefore the maximal invariant set in each bounded set $B$ is necessarily compact.
% hence each bounded subset of $X$ is naturally {\bf admissible} for $\Phi$ in the terminology of Rybakowski \cite{Ryba}.
\vs

A {compact} invariant set $S$ (we allow $S=\emp$) of $\Phi$ is said to be {\em isolated}, if there exists a  neighborhood $N$ of $S$ such that $S$
is the  maximal invariant set in $ \ol N$. Correspondingly, $N$ is called
 an {\em isolating neighborhood} of $S$.

Let there be given an isolated  invariant set $S$.
 A pair of bounded closed subsets $(N,E)$ is called an {\em  index pair} of $S$, if (1) $N\sm E$ is an isolating neighborhood of $S$;
(2) $E$ is $N$-invariant, i.e., for any $x\in E$ and $t\geq 0$, $$\Phi([0,t])x\subset
N\Longrightarrow\Phi([0,t])x\subset E;$$
(3) $E$ is an exit set of $N$. Namely,  for any $x\in N$, if $\Phi(t_1)x\not\in N$ for some
$t_1>0$, then there exists  $0\leq t_0\leq t_1$ such that
    $\Phi(t_0)x\in E.$

%    Finally, let us also recall the concept of an isolating block.
\vs
Let $B\subset X$ be a closed domain. A point $x\in \pa B$  is called a {\em strict egress} (resp. {\em strict ingress}, {\em bounce-off}) point of $B$, if for every trajectory $\gamma:[-\tau,s]\ra X$ with $\gamma(0)=x$ (where $\tau\geq0$, and $s>0$),
\begin{enumerate}
\item[(1)] there exists  $0<\ve<s$ such that
$$
\gamma(t)\not\in B \,\,\,(\mb{resp. }\, \gamma(t)\in \mb{int}B,\,\,\,\mb{resp. }\,\gamma(t)\not\in B),\Hs \A\,t\in (0,\ve);
$$
\item[(2)] if $\tau>0$ then there exists $0<\de<\tau$ such that
$$
\gamma(t)\in \mb{int}B \,\,\,(\mb{resp. }\, \gamma(t)\not\in B,\,\,\,\mb{resp. }\,\gamma(t)\not\in B),\Hs \A\,t\in (-\de, t).
$$
\end{enumerate}
Denote  $B^e$ (resp. $B^i$, $B^b$) the set of all strict egress (resp. strict ingress, bounce-off) points of the closed set $B$, and set 
\be\label{bpm}B^-=B^e\cup B^b,\hs B^+=B^i\cup B^b.\ee
A bounded closed domain $B$ is called an {\em isolating block} \cite{Ryba}, if $B^-$ is closed and $
\pa B=B^i\cup B^-.$ For an isolating block $B$,  we infer from \cite{Ryba} that $(B,B^-)$ is an index pair of the maximal compact invariant set $S$ (possibly empty) in $B$.
%For convenience in statement,  we will call $B^-$ the {\em boundary exit set} of $B$.

%\br %We infer from \cite{Ryba}  that any isolating neighborhood $N$ of $S$contains an index pair.
%An index pair in the terminology of \cite{Ryba} needn't be bounded. However,  bounded index pairs are sufficient for our purposes here.\er

\begin{definition} (homotopy index) Let $(N,E)$ be an index pair of $S$. The
{ homotopy Conley index} of $S$ is defined to be the homotopy
type $[(N/E,[E])]$ of  the pointed space $(N/E,[E])$, denoted by $h(\Phi,S)$.
\end{definition}
\br
Denote  $H_*$ and $H^*$ the singular homology and cohomology theories with coefficient group $\mathbb Z$, respectively.
Applying $H_*$ and $H^*$ to $h(\Phi,S)$ one immediately obtains the { homology} and { cohomology Conley index} $CH_*(\Phi,S)$ and $CH^*(\Phi,S)$, respectively.
\er

An important property of the Conley index is the homotopy invariance  of the index. Here we state a result in this line for the reader's convenience, which is actually a particular case of \cite[Chap.\,I,  Theorem 12.2]{Ryba}.

Let $\Phi_\lam$ ($\lam\in \Lam$) be a family of
asymptotically compact local semiflows on $X$, where $\Lam$ is a metric space.
We say that $\Phi_\lam$ {\it depends on $\lam$
continuously}, if $\Phi_\lam(t)x$ is defined at $(t,x,\lam)$, then  for any sequence $(t_n,x_n,\lam_n)$ converging to
 $(t,x,\lam)$, $\Phi_{\lam_n}(t_n)x_n$ is also defined for all $n$ sufficiently
large, furthermore,  $$\Phi_{\lam_n}(t_n)x_n\rightarrow \Phi_\lam(t)x\hs \mb{ as }\,n\rightarrow\infty.$$

Suppose the family $\Phi_\lam$ ($\lam\in \Lam$) depends on $\lam$ continuously.
Set
$$
  \Pi(t)(x,\lambda)=(\Phi_\lambda(t)x,\lambda), \Hs (x,\lambda)\in
 \sX:= X\times\Lam,\,\,t\geq0.
  $$
Then $\Pi$ is a local semiflow on the product space $\sX$. For convenience, we call $\Pi$  the {\it skew-product flow}  of the family.

We say that  $\Phi_\lam$ is {\em $\lam$-locally uniformly asymptotically compact} ($\lam$-l.u.a.c. in short), if the skew-product flow $\Pi$ is asymptotically compact. It is easy to see
  that if  $\Phi_\lam$ is $\lam$-l.u.a.c., then the parameter space  $\Lam$ is necessarily { locally compact}.

 Suppose for each $\lam\in\Lam$, $\Phi_\lam$ has an isolated invariant set $S_\lam$.
 We call the pair $(\Phi_\lam,S_\lam)$ a {\em dynamic continuation} on $\Lam$, if  for every $\lam_0\in\Lam$, there is a neighborhood $W$ of $\lam_0$ in $\Lam$ and a set $N\subset X$ such that
 $N$ is an isolating neighborhood of $S_\lam$ for all $\lam\in W$.
In the case where $\Phi_\lam$ is $\lam$-l.u.a.c., one can easily verify that a dynamic continuation $(\Phi_\lam,S_\lam)$ is {$\cS$-continuous} in the terminology of \cite[Chap.\,I,  Def. 12.1]{Ryba}. Hence by  \cite[Chap.\,I,  Theorem 12.2]{Ryba} we have

\bt\label{t2.8} \,Assume that $\Phi_\lam$ is $\lam$-l.u.a.c. Let  $(\Phi_\lam,S_\lam)$ be a {dynamic continuation} on $\Lam$. Then  $h(\Phi_\lam,S_\lam)$ is constant for $\lam$  in any component of $\Lam$.
\et
%\noindent{\bf Proof.} Since $\Phi_\lam$ ($\lam\in\Lam$) is $\lam$-l.u.a.c., one can easily verify that the pair $(\Phi_\lam,S_\lam)$  is {$\cS$-continuous} in the terminology of \cite{Ryba}. The conclusion then immediately follows from \cite[Chap.\,I,  Theorem 12.2]{Ryba}. $\Box$

%\br\label{r2.7}Assume that $\Phi_\lam$  is $\lam$-l.u.a.c. Let  $(\Phi_\lam,S_\lam)$ be a {dynamic continuation} on $\Lam$. Then by the $\lam$-l.u.a.c. property of $\Phi_\lam$ and \er

\section{Equilibrium Index of Statically Isolated  Sets}
Let $X$ be a Banach space with norm $\|\.\|$, and $A$ be a sectorial operator on $X$ with {compact} resolvent. Denote  $X^\b$ ($\b\in\R$) the fractional powers of $X$ induced by $A$ equipped with the usual norm $\|\.\|_\a$ (see \cite[Chap. 1]{Henry} for details).
Consider the equation
 \be\label{e3.4}
u_t+Au=f(u), \Hs u=u(t)\in U,
\ee
where $U$ is an open subset of $X^\a$ for some fixed $\a \in[0,1)$. Our main purpose in this section is to introduce the notion of {\em equilibrium index} for \ef{e3.4}.

For convenience in statement, we will write $B\ss U$, meaning that $B$ is a {\bf bounded} subset of $U$ with
$$d(B,\pa U)>0.$$

We   always assume $f:U\ra X$ is locally Lipschitz continuous, namely,
\benu
\item[{\bf (LC)}]   for any $B\ss U$, there exists  $L>0$ such that
$$
||f(x)-f(y)||\leq L||x-y||_{\a},\Hs \A\,x,y\in B.
$$
\eenu
 Under the above assumption, it is well known that the initial value problem of the equation is well-posed. That is,  for each $u_0\in X^\a$ the problem has a unique  solution $u(t)=u(t;u_0)$ with $u(0)=u_0$ defined on a  maximal existence interval $[0,T_{u_0})$ with $u(0)=u_0$; see e.g. \cite[Chap.\,3]{Henry}.
Define a local semiflow on $U$ as
$$\Phi(t)u_0=u(t;u_0),\Hs u_0\in U,\,\,t\in[0,T_{u_0}).
$$ $\Phi$ is usually called the local semiflow generated by \eqref{e3.4}.

\br Note that in general the domain  $U$ may not  be complete. However, one easily understands that the dynamical systems theory reviewed in Section 2 applies to  $\Phi$ quite well on any closed domain $N\ss U$.
\er

\subsection{Definition of the equilibrium index}

Given $K\subset U$, denote $\sE_\Phi(K)$ the set of equilibrium points of $\Phi$ in $K$.
\begin{definition} A set $S\ss U$ is called statically isolated if it has  a neighborhood $N\ss U$ such that
$
\sE_\Phi(\ol N)=\sE_\Phi(S).
$
Correspondingly   $N$  is called a statically isolating neighborhood of $S$.
\end{definition}

%\br If $S\subset U$ is a statically isolated invariant set, then by compactness of $S$ the set $N $ in the above definition can be chosen to be bounded with $\ol N\subset U$.\er

For convenience  in statement, given $Z\subset \mathbb{C}$ and $a\in\R$, we write $$\mb{Re}\,Z\leq a\,\,(resp.\,\,\geq a),$$   meaning that
$
\mb{Re}\,z\leq a\,\,(resp.\,\,\geq a)$ for all $z\in Z.$

Take a number $a\in\R$ with $\mb{Re}\,\sig(aI+A)\geq\de>0$. Let
$$
F:=(aI+A)^{-1}(aI+f).
$$
Then $F(B)$ is pre-compact in $X^\a$ for any $B\ss U$.

\begin{definition} Let $S\ss U$ be statically  isolated with $N\ss U$ being an isolating neighborhood. The equilibrium  index of $S$, denoted by $\mb{\em Ind}\,(\Phi,S)$, is defined as
$$
\mb{\em Ind}\,(\Phi,S)=\deg\,(I-F,\ol N,0),
$$
where $\deg\,(\.,\.,\.)$ denotes the Leray-Schauder degree.
\end{definition}

It is trivial to verify that  the definition of $\mb{Ind}\,(\Phi,S)$ is independent of the choice of $N$ and the number $a$.

A basic property of the equilibrium index  is that it is invariant under transformations.
Specifically, let $Y$ be another  Banach space. Suppose there is an isomorphism $T:X\ra Y$ (a bounded linear bijection with bounded  inverse).
Consider the following equation:
\be\lb{et}
v_t+Bv=g(v),\Hs v=v(t)\in V=TU,
\ee  where $B=TAT^{-1}$, and $g=TfT^{-1}.$
Denote $\Psi$ the local semiflow generated by \ef{et}. Then $\Psi$ and $\Phi$ are conjugate, namely, $\Psi T=T\Phi$. We have
\bp\lb{p3.4}Let $S\ss U$ be a statically isolated  set of $\Phi$. Then $TS\ss V$ is a statically  isolated set of $\Psi$; furthermore,
$$
\mb{\em Ind}\,(\Phi,S)=\mb{\em Ind}\,(\Psi,TS).
$$
\ep
{\bf Proof.} Take a number $a>0$ so that $\mb{Re}\,\sig(A_a)\geq\de>0$. Let
$f_a=aI+f.$ Then \ef{et} can be written  as
\be\lb{etb}
v_t+B_av=g_a(v),
\ee
where
$$
B_a=aI+B=TA_aT^{-1},\hs g_a=aI+g=T f_aT^{-1}.
$$
We have
$$
\mb{Re}\,\sig(B_a)=\mb{Re}\,\sig(A_a)\geq\de>0.
$$

Pick an (statically) isolating neighborhood $N$ of $S$. Then $M=TN\ss V$ is an isolating neighborhood of $TS$. Hence by definition we have
$$\ba{ll}\mb{Ind}\,(\Psi,TS)&=\mb{deg}\,(I-B_a^{-1}g_a,\ol M,0)\\[1ex]&=\mb{deg}\,\(I-T(A_a^{-1}f_a)T^{-1},\ol M,0\)\\[1ex]
&=\mb{deg}\,\(I-A_a^{-1}f_a,\ol N,0\)=\mb{Ind}\,(\Phi,S),\ea
$$
which completes the proof of the result. $\Box$

\subsection{Continuation property}\label{s3.2}
%One of the most important properties of the equilibrium index is the {\em continuation property}. To give a precise description,
Let $X$, $U$ and $A$ be as above. Consider the equation with parameter $\lam\in\Lam$:
\be\label{e3.2}
u_t+A u=f_\lam(u), \Hs u=u(t)\in U,
\ee
  where $\Lam$ is  a {\em connected  compact} metric space. Instead of (LC) we assume that  %$f_\lam(x)$ satisfies  the following local Lipschitz continuity condition in $(x,\lam)$:
 \benu
\item[{\bf (UL)}]   for any   $B\ss U$ and compact subset $\Lam_0$ of $\Lam$, there exists  $L>0$ such that
$$
||f_\lam(x)-f_{\lam'}(y)||\leq L\(||x-y||_{\a}+d(\lam,\lam')\),\Hs (x,\lam),\,(y,\lam')\in B\X\Lam_0.
$$
%the function $f_\lam(x)$ is Lipschitz in $(x,\lam)$ on $B\X\Lam_0$.
\eenu
Denote $\Phi_\lam$ the local semiflow  generated  by \ef{e3.2} on $U$.
It is easy to  verify that $\Phi_\lam$ depends on $\lam$ continuously.
 %Also, by very standard argument (see e.g.\,\,\cite[Chap.\,I, Theorem 4.4]{Ryba}), it can be shown that   %$\{\Phi_\lam\}_{\lam\in\Lam}
  %$\Phi_\lam$ is $\lam$-l.u.a.c.

  Let $\Pi$ be skew-product flow of the family $\Phi_\lam$ ($\lam\in\Lam$) on $\cU=U\X\Lam$. Clearly $x$ is an equilibrium of $\Phi_\lam$ {iff}\, $(x,\lam)$ is an equilibrium of $\Pi$.

%Suppose  $\Phi_\lam$ has a statically isolated   set $S_\lam\ss U$ for each $\lam\in\Lam$ (we allow $S_\lam=\emp$).

%\begin{definition} We say that $S_\lam$ is a static continuation of the family $\Phi_\lam$  $(\lam\in\Lam)$ on $\Lam$, if\,(1) \,$S_\lam$ is upper semicontinuous in $\lam$; and (2) \,
%  for each $\lam_0\in \Lam$, there exists $N\ss U$  such that $N$ is a statically isolating neighborhood of $S_\lam$ for all $\lam$ in a neighborhood $W$ of $\lam_0$.
%  \end{definition}
The following result is a simple consequence of the continuation property of the Leray-Schauder degree.

  \bt\label{t3.5} Let $\cS$ be a statically isolated set of the skew-product flow $\Pi$ in $\cU$. Set   $S_\lam:=\cS[\lam]$.  Then $\mb{\em Ind}\,(\Phi_\lam,S_\lam)$ remains constant for $\lam\in \Lam$.
\et

  \noindent{\bf Proof.}  Since $\cS$ is a statically isolated set of $\Pi$, by definition $\cS$ has a bounded  closed isolating neighborhood $\cN\Subset \cU$. Note that $d(\cN,\pa_\sX\cU):=d>0$, where $\pa_\sX\cU$ denotes the boundary of $\cU$ in $\sX:=X^\a\X\Lam$. It is trivial to check that $d(\cN,\pa_\sX\cU)\leq d(\cN[\lam],\pa U)$. Therefore we deduce that
  $$
  d(\cN[\lam],\pa U)\geq d>0,\Hs \lam\in\Lam.
  $$
 Hence we can  pick a bounded closed set $B\Subset U$ such that $\cN[\lam]\subset B$ for all $\lam\in \Lam$.

 Using  local Lipschitz continuity of $f_{\lam}$  one   easily   verifies that
$\ba{ll}\Cup_{\lam\in \Lam}f_\lam\(B\)\ea$ is  bounded in $X$.
 Consequently $M:=\Cup_{\lam\in \Lam}(aI+f_\lam)(B)$ is a bounded set in $X$.
 Let $$F_\lam=(aI+A)^{-1}(aI+f_\lam).$$ Then the set
$\ba{ll}\Cup_{\lam\in \Lam}F_\lam(B)=(aI+A)^{-1}M \ea
  $
  %\be\label{e3.6}F_\lam=(aI+A)^{-1}(aI+f_\lam).\ee
is bounded in $D(A)=X^1$. As  $A$ has compact resolvent, the embedding $X^1\hookrightarrow X^\a$ is compact. Therefore $\Cup_{\lam\in \Lam}F_\lam(B)$ is precompact in $X^\a$.
Now  by  virtue of continuation  property of the Leray-Schauder degree we immediately concludes that
$$\mb{Ind}(\Phi_\lam,S_\lam)=\mb{deg}\,(I-F_\lam,\cN[\lam],0)\equiv \mb{const.}$$
 on $\Lam$. The proof of the theorem is complete. $\Box$
 \Vs

 Suppose  $\Phi_\lam$ has a statically isolated invariant  set $S_\lam$ for each $\lam\in\Lam$ (we allow $S_\lam=\emp$).
 We say that $(\Phi_\lam,S_\lam)$ is a {\em static continuation} on $\Lam$, if for each $\lam_0\in \Lam$, there exist $N\Subset U$ and  a neighborhood $W$ of $\lam_0$ in $\Lam$ such that $N$ is a statically isolating neighborhood of $S_\lam$ for all $\lam\in W$.

\bt\label{t3.8}
Let  $(\Phi_\lam,S_\lam)$ be a static continuation on $\Lam$.
Then $\mb{\em Ind}\,(\Phi_\lam,S_\lam)$ remains constant on  $\Lam$.
In particular, if $\Lam=[0,1]$ then
$
\mb{\em Ind}\,(\Phi_0,S_0)=\mb{\em Ind}\,(\Phi_1,S_1).
$
\et

\noindent{\bf Proof.} To prove Theorem \ref{t3.8}, it suffices to check that for any $\lam_0\in\Lam$, there is a neighborhood $\cO$ of $\lam_0$ in $\Lam$ such that
\be\label{cpei}\mb{Ind}(\Phi_\lam,S_\lam)\equiv const.,\Hs \lam\in \cO.\ee

Let $\lam_0\in\Lam$. Then by the definition of a static continuation, one can find a neighborhood $N\Subset U$ of $S_{\lam_0}$ and  a neighborhood $W$ of $\lam_0$ in $\Lam$ such that $N$ is a statically isolating neighborhood of $S_\lam$ for all $\lam\in W$. Pick a connected compact neighborhood $\cO$ of $\lam_0$ with $\cO\subset W$. One easily sees that $\cS:=\Cup_{\lam\in\cO}S_\lam\X\{\lam\}$ is a static isolated set of the skew-product flow of $\Phi_\lam$ ($\lam\in\cO$). Applying Theorem \ref{t3.5} we immediately conclude the validity of \ef{cpei}. $\Box$

\subsection{The finite dimensional case}

Now we assume $X$ is an $m$-dimensional Banach space. Consider the ODE system:
\be\label{e3.2}
x'(t)=f(x),\Hs x=x(t)\in U,
\ee
where $f:U\ra X$ is locally Lipschitz. Denote $\Phi$ the  semiflow generated by   \eqref{e3.2}.

\bp\label{p3.3} Let $S$  be a statically isolated set of $\Phi$. Then
$$
\mb{\em Ind}\,(\Phi,S)=(-1)^m \deg(f,N,0)
$$
for any  isolating neighborhood $N$ of $S$.

If $S$ is an isolated invariant set of $\Phi$, we also have
$$
\mb{\em Ind}\,(\Phi,S)=\chi(h(\Phi,S)),
$$
where $\chi(h(\Phi,S))=\Sig_{q=0}^\8\,(-1)^q\b_q(h(\Phi,S))$ is the Euler number of the Conley index $h(\Phi,S)$.
\ep
{\bf Proof.} We may assume  $X=\R^m$. Take an arbitrary $m\X m$-matrix  $A$ and rewrite \eqref{e3.2} as
$$
x'(t)+Ax=\~f(x),
$$
where $\~f(x)=Ax+f(x)$. Pick a number $a>0$ so that $\mb{Re}\,\sig(aI+A)>0$. Then
$$\ba{ll}
\mb{Ind}(\Phi,S)&=\mb{deg}\(I-(aI+A)^{-1}(aI+\~f),\,N,0\)\\[1ex]
&=\mb{deg}\(I-(aI+A)^{-1}(aI+A+f),\,N,0\)\\[1ex]
&=\mb{deg}\(-(aI+A)^{-1}f,\,N,0\).
\ea
$$
Note that $\mb{Re}\,\sig(aI+A)>0$ implies $\mb{det}\,(aI+A)>0$. Further by the definition of the Brouwer degree it is easy to deduce that
 $$
\mb{deg}\(-(aI+A)^{-1}f,\,N,0\)=\mb{deg}\(-f,\,N,0\)=(-1)^m \mb{deg}\(f,\,N,0\).
$$
Hence
$$
\mb{Ind}(\Phi,S)=(-1)^m \mb{deg}\(f,\,N,0\).
$$

Now assume that  $S$ is an  isolated invariant set of $\Phi$. Pick an (dynamically) isolating neighborhood $N$ of $S$. Then we  infer from Reineck \cite{Rein} (see also \cite[Theorem 3.2]{Raz}) that $S$ can be continued in $N$ to an isolated invariant set $K$ of a flow $G$  generated by
$$
x'(t)=g(x),
$$
where $g$ is a Morse-Smale gradient vector field on a neighborhood $V$ of $K$ with  $V\subset N$, and $g(x)=f(x)$ on $\R^m\sm N$. Thus by Rybakowski \cite[Chap. III,Theorem 3.8]{Ryba} (see also Chang \cite[Chap. II, Theorems 3.1-3.3]{Ch1}) one concludes that
$$\ba{ll}
\mb{Ind}(\Phi,S)&=(-1)^m \mb{deg}\(f,\,N,0\)=(-1)^m \mb{deg}\(g,\,N,0\)\\[1ex]
&=\Sig_{q=0}^\8(-1)^q\b_q(h(G,K))\\[1ex]
&=\chi(h(G,K))=\chi(h(\Phi,S)).\ea
$$
The proof is complete. $\Box$
%where $\b_q(h(\Phi,S))$ is the Betti number of $h(\Phi,S)$.

\br In \cite{Raz} the Euler number $\chi(h(\Phi,S))$ is defined to be the Poincar\'{e} index of $S$ (see \cite[Section 3]{Raz}).
\er
\section{ Invariant Manifolds of Nonlinear Equations}
For completeness and the reader's convenience, in this section we briefly recall some results on local invariant manifolds of the following nonlinear evolution equation with parameter $\lam\in\R$:
\be\label{e4.1}
u_t+A u=f_\lam(u),
\ee
where $A$ is a sectorial operator in $X$ with compact resolvent, and $f_\lam(x)$ is a continuous  mapping from $X^\a\X \R$ to $X$ for some $\a\in[0,1)$.
% which is differentiable  in $u$ with $f_\lam'(u)$ being continuous in $(u,\lam)$.
%As in Section 3, we always assume that the local Lipschitz continuity condition (LC) is fulfilled by $f=f_\lam(u)$ for all  $\lam\in\R$.

Without loss of generality, assume that
$
f_\lam(0)\equiv0$
for $\lam\in\R,$ hence $u=0$ is always an equilibrium solution of  (\ref{e4.1}).

Assume that $f_\lam(x)$ is locally Lipschitz in $x$; furthermore, the local Lipschitz continuity assumption  (UL) in Subsection \ref{s3.2}  is fulfilled. Denote  $\Phi_\lam$ the local semiflow generated by \ef{e4.1}. Then we know that $\Phi_\lam$ depends on $\lam$ continuously.
 Also, by very standard argument (see e.g.\,\,\cite[Chap.\,I, Theorem 4.4]{Ryba}),  it can be shown that $\Phi_\lam$ is $\lam$-l.u.a.c. on $X^\a$.

\subsection{Fundamental assumptions and notations}
Let $L_{\lam}=A-f_\lam'(0)$, and write $$g_\lam(u)=f_\lam(u)-f_\lam'(0)u.$$  Then  \ef{e4.1} reads
\be\lb{e3.1}
u_t+L_\lam u=g_\lam(u).
\ee
%Since $g_\lam:E\ra X$ is locally Lipschitz, the initial value problem of (\ref{e3.1}) is well-posed \cite{Henry}. That is, for each  $u_0\in E$ the problem has a unique solution $u(t)$ in $E$ with $u(0)=u_0$ on some maximal existence interval $[0,T_0)$.
Suppose  there exists $\nu>0$ such that  the following hypotheses (H1)-(H3) are fulfilled for every $\lam\in J_0=[\lam_0-\nu,\,\lam_0+\nu]$:
\begin{enumerate}
\item[(H1)] The   spectral $\sig(L_\lam)$  has a decomposition $\sig(L_\lam)=\sig_\lam^1\cup \sig_\lam^2\cup \sig_\lam^3$
%\,, where $\sig_\lam^1\cup \sig_\lam^2$ is finite,
    such that
 \be\label{e:2.4}
\mb{Re}\,\sig_\lam^1\leq -2\de,\hs -\de\leq \mb{Re}\,\sig_\lam^2 \leq \de,\hs \mb{Re}\,\sig_\lam^3\geq 2\de\ee for some  $\de>0$ (independent of $\lam\in J_0$).

\item[(H2)] The space $X$ has a   decomposition
$X=X^1_\lam\oplus X^2_\lam\oplus X^3_\lam\,,$ corresponding to the decomposition of $\sig(L_\lam)$ in (H1), such that $X^i_\lam$  $(i=1,2,3)$ are $L_{\lam}$-invariant subspaces of $X$. Moreover,
\be\label{fie3}
\mb{dim}\,(X^i_\lam)\equiv\mb{dim}\,(X^i_{\lam_0}):= \mathfrak{m}_i<\8,\Hs i=1,2.
\ee

\item[(H3)] The projection operators $$P_\lam^{i}:X\ra X^i_\lam,\hs i=1,2$$ are continuous in $\lam$.

\end{enumerate}\Vs

\br\label{r4.1}  The above assumptions implies that there is  a  family $T_\lam$ ($\lam\in J_0$) of isomorphisms  on $X$  depending continuously on $\lam$ with
$T_{\lam_0}=I$, such that \be\label{IS}
    T_\lam X^i_\lam=X^i_{\lam_0}:=X^i,\Hs i=1,2,3;
\ee
see Appendix A for the proof. It is trivial to verify that
$
T_\lam P_\lam^i=P_{\lam_0}^i$ for all $\lam\in J_0.$
\er

Let $$X_\lam^{ij}=X_\lam^i\oplus X_\lam^j,\Hs 1\leq i,j\leq 3,\,\,\,i\ne j,$$
 and denote  $P_\lam^{ij}$ the projection from $X$ to $X_\lam^{ij}$.
 Then $P_\lam^{ij}=P_\lam^{i}+P_\lam^{j}$. We infer from  Remark \ref{r4.1} that
 \be\label{e:4.35}
T_\lam P_\lam^{ij}=P_{\lam_0}^{ij},\Hs \A\,\lam\in J_0.
\ee

\vs  We will  rewrite  $X^\a=E$. Let
$$
E^i_\lam=E\cap X^i_\lam,\hs E^{ij}_\lam=E\cap X^{ij}_\lam\,\,\,(i\ne j).
$$
Then
$$E=E^1_\lam\oplus E^2_\lam\oplus E^3_\lam,\hs\, E^{ij}_\lam=E^i_\lam\oplus E^j_\lam.$$
 Because $X^1_\lam$ and $X^2_\lam$ are finite dimensional, we actually  have $$E^1_\lam=X^1_\lam,\hs E^2_\lam=X^2_\lam,\hs E^{12}_\lam=X^{12}_\lam.$$

For notational  simplicity,  hereafter  we write $$X_{\lam_0}^i=X^i,\hs  X_{\lam_0}^{ij}=X^{ij},\hs  E_{\lam_0}^i=E^i,\hs  E_{\lam_0}^{ij}=E^{ij},$$
$$P^i=P_{\lam_0}^i,\hs P^{ij}= P_{\lam_0}^{ij}.$$%  the projection operators from $E$ to $E^i$.
\vs
 Denote   $\mB^i_r$ the ball in  $E^i$  centered at $0$ with radius $r$,
%$$\mB_r=\{x\in E:\,\,||x||_\a<r\} ,\hs \mB^i_r=\{x\in E^i:\,\,||x||_\a<r\},$$
and  set
$$
  \Xi(r)=\mB^1_r\oplus \mB^2_r\oplus \mB^3_r.
  $$
%The following easy fact will be frequently used.
  \bl\label{p4.1} For any neighborhood $\cU$ of $0$ in $E$, there exists $r>0$ such that
\be\label{rt0}
\Xi(r)\subset \cU.
\ee
\el
{\bf Proof.}
 Since $P^i\,\cU$ is a neighborhood of $0$ in $E^i$, we can  pick an  $\ve>0$ sufficiently small  so that $\ol\mB^i_\ve\subset P^i\cU$ for $i=1,2$.  Further  by using compactness of $\ol\mB^i_\ve$ ($i=1,2$) one can easily verify   that there exists  $\ve'>0$ such that
$\ol\mB^1_\ve\oplus \ol\mB^2_\ve\oplus \ol\mB^3_{\ve'}\subset \cU$, from which
 \ef{rt0} immediately follows. $\Box$

\subsection{Local invariant manifolds}

Let us recall briefly  some fundamental results on
local invariant manifolds.

Let  $T_\lam$ be the isomorphism on $X$ given  in Remark \ref{r4.1}. Setting $v=T_\lam u$,
system \eqref{e3.1} is transformed  into an equivalent one:
\be\label{rt2} v_t+ B_\lam v=h_\lam(v),\Hs v=v(t)\in E=X^\a\ee
where \be\label{e3.1a}B_\lam=T_\lam L_\lam T^{-1}_\lam,\hs h_\lam=T_\lam g_\lam T^{-1}_\lam.\ee
 It is easy to verify (or see \cite[Section 3]{LW}) that
$\sig( B_\lam)=\sig(L_\lam).$

Denote $\Psi_\lam=\Psi_\lam(t)$ the local semiflow generated by \ef{rt2}.

Since $X^i_\lam$ ($i=1,2,3$) are $L_{\lam}$-invariant, by \eqref{IS} we find   that  $X^i:=X^i_{\lam_0}$ ($i=1,2,3$) are  $ B_\lam$-invariant for all $\lam\in J_0$. Thus by  (\ref{e:2.4}) we deduce that
\be\label{e:4B}\mb{Re}\,\sig( B_{\lam}^1)\leq -2\de,\hs -\de\leq\mb{Re}\,\sig( B_{\lam}^2) \leq \de,\hs \mb{Re}\,\sig( B_{\lam}^3)\geq 2\de\ee
where $ B^i_\lam= B_\lam|_{X^i}$  is the restriction of $ B_\lam$ on $X^i$. Now the same arguments as in the proofs of Lemmas 3.3 and 3.4 in \cite{LW} apply to prove the following theorem.
%\vs Because $g'_\lam(0)=0$ and $g'_\lam(v)$ is continuous in $(v,\lam)$, it is easy  to  see  that $||h_\lam'(v)||\ra 0$ as $||v||_\a\ra 0$ uniformly with respect to $\lam\in  J_0$. Further by the Mean-value Theorem it can be shown  that for any $\ve>0$, there is a neighborhood $\cO$ of $0$ in $X^\a$ such that
%\be\label{rt3}
%||h_\lam(u)-h_\lam(v)||\leq \ve ||u-v||_\alpha,\Hs \A\,u,v\in \cO,\,\,\lam\in J_0.
%\ee

\bt\label{rtt0} There exist open convex neighborhoods $V_i$ of $0$ (independent of $\lam\in J_0$) in $E^i$\, $(i=1,2,3)$ such that the following assertions hold.
\benu
\item[$(1)$] There is a continuous mapping $\psi_{\lam}(x)$ from $(V_1\oplus V_2)\X  J_0$ to $V_{3}$  which is differentiable in $x$ with $\psi_\lam'(x)$ being continuous in $(x,\lam)$ and  $\psi_\lam'(0)=0$ such that
     for each $\lam\in J_0$,
the set $M^{12}_\lam=\{x+\psi_{\lam}(x):\,\,x\in V_1\oplus V_2\}$
 is a local invariant manifold of  system \eqref{rt2}.
  %Consequently $${\cM}^{12}_\lam:=T^{-1}_\lam M^{12}_\lam$$ is a local invariant manifold of the original system (\ref{e4.1}).

\item[$(2)$] There is a  continuous mapping  $\rho_{\lam}(x_2)$ from $V_2\X  J_0$ to $V_{1}$ which is differentiable in $x_2$ with $\rho_\lam'(x_2)$ being continuous in $(x_2,\lam)$  and $\rho_\lam'(0)=0$ such that for each $\lam\in J_0$, the set
$ M^{2}_\lam=\{x_2+\zeta_\lam(x_2):\,\,x_2\in V_2\}$
     is a local invariant manifold of  \eqref{rt2}, where
      \be\lb{rt8d}
\zeta_\lam(x_2)=\rho_\lam(x_2)+\psi_\lam(x_2+\rho_\lam(x_2)).
\ee
%Consequently $$\cM^{2}_\lam:=T^{-1}_\lam M^{2}_\lam$$ is a local invariant manifold of (\ref{e4.1}).

\item[$(3)$] $S\subset V:=V_1\oplus V_2\oplus V_3$ is an isolated invariant set of the semiflow $\Psi_\lam$ of \eqref{rt2}  {\em iff }  $S\subset M^{12}_\lam$ (\,resp. $M_\lam^2$\,)  and is an isolated invariant set of  $\Psi_\lam^{12}=\Psi_\lam|_{M_\lam^{12}}$ (\,resp. $\Psi_\lam^2=\Psi_\lam|_{M_\lam^2}$\,). Furthermore,
%    Consequently  $S\subset U_\lam=T_\lam^{-1}V$ is an isolated invariant set of the semiflow $\Phi_\lam$ of \eqref{e4.1} {\em iff }  $S\subset \cM^2_\lam$  and is an isolated invariant set of  $\Phi_\lam^{2}=\Phi_\lam|_{\cM_\lam^2}$. Furthermore,
    \be\label{ertc}
    h(\Psi_\lam,S)=h(\Psi_\lam^{12},S)=\Sigma^{\mathfrak{m}_1}\wedge h(\Psi_\lam^{2},S).
    \ee

\eenu
\et
\br Assertions (1) and (2) are only  slight  modifications of some classical results in the geometric theory of evolution equations (see e.g. \cite[ Chap. II, Theorem 2.1]{Ryba}).
A sketch of the proof can be found in \cite{LW}.  (3)  is a parameterized version of a corresponding result in \cite[Chap. II, Theorem 3.1]{Ryba}. \er

%Consequently $\cM^2_\lam=T^{-1}M^2_\lam$ is a local invariant manifold of the system (\ref{e:1.1}).
\br\label{r:4.5}
Using  $\psi_\lam'(0)=0$ and $\rho_\lam'(0)=0$ it is easy to deduce that for any $\ve>0$, there exists $r>0$ such that
$$||\psi_\lam'(x)||_\a<\ve,\hs ||\rho_\lam'(x_2)||_\a<\ve
$$
for all $x\in $\mB$^1_r\oplus\mB^2_r$, $x_2\in\mB^2_r$ and $\lam\in J_0$.
It follows that
\be\label{rt6}
\psi_\lam(x)=o(||x||_\a^2)\,\,\,(||x||_\a\ra0),\hs
\rho_\lam(x_2)=o(||x_2||_\a^2)\,\,\,(||x_2||_\a\ra0)
\ee
uniformly with respect to $\lam\in J_0$.
\er

\br It is known (see e.g. \cite[Section 2.2]{Ryba}) that  $\psi_\lam$ satisfies equation
\be\label{rt8b}\ba{ll}
\psi_\lam'(x)\[ B_\lam^{12}x-P^{12}h_\lam(x+\psi_\lam(x))\]
= B_\lam^3\psi_\lam(x)-P^3h_\lam(x+\psi_\lam(x))\ea
\ee for $x\in V_1\oplus V_2$, and $\rho_\lam$ satisfies
\be\label{rt8c}\ba{ll}
\rho_\lam'(x_2)\[ B_\lam^{2}x_2-P^{2}h_\lam(x_2+\zeta_\lam(x_2))\]= B_\lam^1\rho_\lam(x_2)-P^1h_\lam(x_2+\zeta_\lam(x_2))\ea
\ee
 for $x_2\in V_2$. Here (and below) $P^i=P_{\lam_0}^i,$ $P^{ij}= P_{\lam_0}^{ij},$ and
  $$ B_\lam^{i}= B_\lam|_{X^{i}},\hs B_\lam^{ij}= B_\lam|_{X^{ij}}.$$
    \er

Let $$U_\lam^{12}=T^{-1}_\lam (V_1\oplus  V_2),\hs U_\lam^2=T^{-1}_\lam V_2.$$ Define $\phi_\lam:U_\lam^{12}\ra E^3_\lam$ and    $\xi_\lam:U_\lam^2\ra E^{13}_\lam$ as below:
\be\label{phi}\phi_\lam(x)=(T^{-1}_\lam \psi_\lam T_\lam)(x),\Hs x\in U_\lam^{12},\ee
\be\label{zeta}
\xi_\lam(x_2)=(T^{-1}_\lam \zeta_\lam\, T_\lam)(x_2),\Hs x_2\in U_\lam^2.
\ee
As a direct consequence of Theorem \ref{rtt0}, we have the following.
\bt\label{t:4.7} Assume the hypotheses  (H1)-(H3).  Then
\benu
\item[$(1)$] the following two sets are local invariant manifolds of system \eqref{e4.1}:
$${\cM}^{12}_\lam:=T^{-1}_\lam M^{12}_\lam=\{x+\phi_{\lam}(x):\,\,x\in U_\lam^{12}\},$$
$${\cM}^{2}_\lam:=T^{-1}_\lam M^{2}_\lam=\{x_2+\xi_{\lam}(x_2):\,\,x\in U_\lam^{2}\};$$
 \item[$(2)$]   $S\subset U_\lam=T_\lam^{-1}V$ (\,where $V=V_1\oplus V_2\oplus V_3$\,) is an isolated invariant set of the semiflow $\Phi_\lam$ of \eqref{e4.1} {\em iff }  $S\subset \cM^{12}_\lam$  (\,resp. $S\subset \cM^2_\lam$ \,) and is an isolated invariant set of  $\Phi_\lam^{12}=\Phi_\lam|_{\cM_\lam^{12}}$
     (\,resp. $\Phi_\lam^{2}=\Phi_\lam|_{\cM_\lam^2}$\,).  Furthermore,
    \be\label{erth}
    h(\Phi_\lam,S)=h(\Phi_\lam^{12},S)=\Sigma^{\mathfrak{m}_1}\wedge h(\Phi_\lam^{2},S).
    \ee
    \eenu
\et

Let $u(t)=x(t)+x_3(t)$ be a solution of \ef{e3.1} (or \ef{e4.1}) lying in $\cM_\lam^{12}$, where $$x=x(t)\in U_\lam^{12},\hs x_3=x_3(t)\in U_\lam^3:=T^{-1}_\lam V_3.$$ Then
$x_3=\phi_\lam(x)$, and therefore  $x$ satisfies
\be\label{re12}
\dot{x}+L_\lam ^{12}x=P^{12}_\lam g_\lam\(x+\phi_\lam(x)\),\Hs x=x(t)\in U_\lam^{12}.
\ee
Similarly one can also obtain the equation   corresponding to $\cM_\lam^{2}$:
\be\label{re2}
\dot{x_2}+L_\lam ^{2}x_2=P^{2}_\lam g_\lam\(x_2+\xi_\lam(x_2)\),\Hs x_2=x_2(t)\in U_\lam^{2}.
\ee
 \ef{re12} and \ef{re2} will be referred to as the {\em reduction  equations} of \ef{e3.1}  (or \ef{e4.1}) on $\cM_\lam^{12}$ and $\cM_\lam^{2}$, respectively.

Denote  $\cR_\lam^{12}$ and $\cR_\lam^2$ the local semiflow on $U_\lam^{12}$ and $U_\lam^2$ generated by \ef{re12} and \ef{re2}, respectively. Then $\cR_\lam^{12}$ and $\cR_\lam^2$ conjugate with $\Phi_\lam^{12}$ and $\Phi_\lam^2$, respectively.  More precisely, we have
\be\label{e:4pa}
\cR_\lam^{12}(t)x=P_\lam^{12}\,\Phi_\lam(t)(x+\phi_\lam(x)),\Hs x\in U^{12}_\lam,
\ee
\be\label{e:4p}
\cR_\lam^2(t)x_2=P_\lam^2\,\Phi_\lam(t)(x_2+\xi_\lam(x_2)),\Hs x_2\in U^2_\lam.
\ee
%Therefore $\cR_\lam^2$ conjugates with the restriction $\Phi_\lam|_{\cM_\lam^2}$ of $\Phi_\lam$ on $\cM_\lam^2$. Due to this reason, we will identify $\cR_\lam^2$  with $\Phi_\lam|_{\cM_\lam^2}$, regardless of the conjugacy between them.
 \br\label{p:4.1}%Let $U_\lam:=T_\lam^{-1}V$ be the neighborhood of $0$ in $E$ given  in Theorem \ref{rtt0} (3). Then
By \ef{e:4p}, \ef{e:4p} and  Theorem \ref{t:4.7} it is easy to see that $S\subset U_\lam$ is an (isolated) invariant set of $\Phi_\lam$ {\em iff} $S_{12}=P_\lam^{12}S$ (resp. $S_2=P_\lam^2S$) is an (isolated) invariant set of $\cR_\lam^{12}$ (resp. $\cR_\lam^2$).\er

\section{Reduction Theorem of Equilibrium  Index near Equilibrium Points}

 Our main purpose in this section is to  establish  a  reduction theorem  for  equilibrium indices of isolated invariant sets of system \ef{e4.1} near equilibrium points.

We follow the same notations as in Section 4 and assume the function $f_\lam(x)$ in \ef{e4.1} satisfies all the regularity hypotheses  in Section 4. Suppose also that
$f_\lam(0)\equiv0$ for $\lam\in\R$. The main result is the following theorem.

\bt\label{rtt1}  Assume the hypotheses  (H1)-(H3) in Section 4 are fulfilled.  Then there exists an open neighborhood $U$ of $0$ in $E$ (independent of $\lam\in J_0$)  such that for every  isolated invariant set $S$ of $\Phi_\lam$ $(\lam\in J_0)$ in $U$, we have
\be\label{ert}\ba{ll}
\mb{\em Ind}\,(\Phi_\lam,S)=\mb{\em Ind}\,(\cR_\lam^{12},P_\lam^{12}S)=(-1)^{\mathfrak{m}_1}\mb{\em Ind}\,(\cR_\lam^2,P^2_\lam S),\ea
\ee
where $\mathfrak{m}_i=\dim\,(X^i)$ $(i=1,2)$.
\et

\noindent
{\bf Proof.} Some  basic ideas and techniques used here are borrowed from \cite{Ryba} (see \cite[Chap. II, Theorem 3.1]{Ryba} and its proof).

  For notational simplicity,  in the following we will  drop the subscript ``$\lam$'' and simply rewrite
 $ B_\lam$, $h_\lam$, $\psi_\lam$, $\rho_\lam$, $\zeta_\lam$ and  $\Psi_\lam$
  as $ B$, $h$, $\cdots$, respectively, unless we need to emphasize  the dependence on $\lam$.
We  split the  argument into several steps.
\vs
{\em Step }1.  We begin with system \ef{rt2}.
Let $V_i$ ($i=1,2,3$) be the neighborhood of $0$ in $E^i$ given in Theorem \ref{rtt0}, and let
$$V=(V_1\oplus V_2) \oplus V_3:=\Om\oplus V_3.$$ For each fixed $\lam\in J_0$, define a family of mappings  $Q_\theta: V\ra E$ ($\theta\in[0,1]$) as
\be\label{e:Q}
Q_\theta(u)=x+(x_3-\theta\psi(x)),\Hs \A \,u=x+x_3\in V=\Om\oplus V_3.
\ee
Note that if $u=x+x_3\in M^{12}=M^{12}_\lam$ then
\be\label{e:Q1}
Q_1(u)=x+(x_3-\psi(x))=x=P^{12}u.
\ee

 Using the same argument as in Step 1 in  the proof of \cite[Chap. II, Theorem 3.1]{Ryba},  one can easily verify that $Q_\theta (V):=\~V_\theta$ is open; furthermore, $Q_\theta$ has a continuous  inverse $Q_\theta^{-1}:\~V_\theta\ra V$ given by
$$
Q_\theta^{-1}(v)=x+(y_3+\theta\psi(x)),\Hs \A\,v=x+y_3\in \~V_\theta,
$$
where $x=P^{12}v\in\Om$, and $y_3=P^3v\in E^3$. (Here and below $P^{i}$ and $P^{ij}$ denote $P^{i}_{\lam_0}$ and $P^{ij}_{\lam_0}$, respectively, except otherwise statement.)
  %Observing that
%$$
%||(Q_\theta-I)(x)||_\a=\theta||\psi(x)||_a=\theta||\psi_\lam(x)||_a=0(||x||_\a^2)
%$$
%$\mb{as }||x||_\a\ra0$ uniformly with respect to $\lam\in J_0$ and $\theta\in[0,1]$, one easily sees that  $V$ can also be chosen sufficiently small in advance so that

Note that $Q_\theta$ depends upon $\theta$ and $\lam$ continuously. Consequently  $\~V_\theta=Q_\theta (V)$ is continuous in ($\theta$, $\lam$) in the sense of Hausdorff distance. Using this simple fact one  easily checks that there is  a neighborhood $\~\cO$ of $0$ in $E$  such that $$\~\cO\subset \~V_\theta,\Hs \A\,\theta\in[0,1],\,\,\lam\in J_0.$$

By Lemma\,\ref{p4.1} we can pick a  number  $r>0$ such that $$\Xi(r):=\mB^1_r\oplus \mB^2_r \oplus \mB^3_r\subset \~\cO.$$ ($\mB^i_r$ denotes the ball in $E^i$).  In view of \ef{rt6}, it can be assumed that $r$ is chosen sufficiently small so that for all $\lam\in J_0$, we have
\be\label{e:4.33}
\psi(x)=\psi_\lam(x)\in \mB^3_r,\Hs \A\,x\in \mB^1_r\oplus \mB^2_r.
\ee

\vs
 Define a local semiflow $\~\Psi_\theta$ on $\~V_\theta$ to be the ``\,image'' of the local semiflow  $\Psi=\Psi_\lam$ of  \eqref{rt2} under $Q_\theta$, namely,
$$
\~\Psi_\theta(t) Q_\theta(u)=Q_\theta \Psi (t)(u),\Hs \A\,u\in V.
$$
Then by \ef{e:4.33} $\~\Psi_\theta$ is well defined on the domain $\Xi(r)$ for all $\theta\in[0,1]$ and $\lam\in J_0$.
Making us of the relation in \ef{rt8b}, it can be shown by some simple computations  that $\~\Psi_\theta$ is precisely the local semiflow generated   by system
\be\label{rt7}\left\{\ba{lll}
\dot{x}+ B^{12} x=N_\theta(x+y_3),\\[1ex]
\dot{y}_3+ B^3 y_3=M_\theta(x+y_3),\ea\right.\Hs v=x+y_3\in \~V_\theta,
\ee
where $x=P^{12}v$, $y_3=P^3v$, $N_\theta=P^{12}h  Q_\theta^{-1}$, and
$$\ba{ll}
M_\theta(x+y_3)=&P^3\left[\,h  Q_\theta^{-1}(x+y_3)-\theta h (x+\psi (x))\, \right]\\[1ex]
&+\,\theta \psi' (x)P^{12}\left[\,h (x+\psi(x))-h  Q_\theta^{-1}(x+y_3)\,\right].\ea
$$
Note that  \eqref{rt7} reduces to  system \eqref{rt2} when  $\theta=0$. Hence $\~\Psi_0$ coincides with $\Psi$.
\vs

Take a  number $\ve>0$ (independent of $\lam$ and $\theta$)  such that $\Xi(\ve)\subset V$  and \be\label{rt25}Q_\theta(\Xi(\ve))\subset \Xi(r),\Hs\A\,\theta\in[0,1],\,\,\lam\in J_0.\ee
 Let $S\subset {\Xi(\ve/2)}$ be an isolated invariant set of $\Psi$, and write
 $$S_\theta=Q_\theta S.$$
  Then \ef{rt25} implies that $S_\theta$ is an isolated invariant set of $\~\Psi_\theta$ in $ \Xi(r)$ for all $\theta\in[0,1]$ (and $\lam\in J_0$).
It is trivial to verify  that $(\~\Psi_\theta,S_\theta)$ is a static coninuation. Thus by Theorem\,\ref{t3.8} we deduce  that
\be\label{rt10}
\mb{Ind}(\Psi,S)=\mb{Ind}(\~\Psi_0,S_0)=\mb{Ind}(\~\Psi_1,S_1).
\ee
\vs{\em Step 2.}
For $\theta=1$,  system \eqref{rt7} reads
\be\label{rt12}\left\{\ba{lll}
\dot{x}+ B^{12} x=N_1(x+y_3),\\[1ex]
\dot{y}_3+ B^3 y_3=M_1(x+y_3),\ea\right.
\ee
where
$$\ba{ll}
M_1(x+y_3)=&P^3\left[\,h  Q^{-1}_1(x+y_3)- h (x+\psi (x))\, \right]\\[1ex]
&+\, \psi' (x)P^{12}\left[\,h (x+\psi(x))-h  Q_1^{-1}(x+y_3)\,\right].\ea
$$
Since $Q_1^{-1}(x+y_3)=x+(y_3+\psi(x))$, it is clear  that
\be\label{rt11}M_1(x+y_3)=0\hs \mb{whenever }\,y_3=0.\ee

Consider on the domain $\Xi(r)$ the homotopy of system  \eqref{rt12}:
\be\label{rt13}\left\{\ba{lll}
\dot{x}+ B^{12} x=N_1(x+\theta y_3),\\[1ex]
\dot{y}_3+ B^3 y_3=\theta M_1(x+y_3),\ea\right.\Hs \theta\in[0,1].
\ee

 \bl\label{rtl1} There exist a small open neighborhood $\~U=\~U_1\oplus \~U_2\oplus \~U_3$ of $0$ in $E$ with $\~U\subset \Xi(r)$ and constants $C,\mu>0$  such that for all $\theta\in [0,1]$ and $\lam\in J_0$, if $u(t)=x(t)+y_3(t)$ is a solution of \eqref{rt13} on $[t_0,t_1]$ lying in $\~U$ then
$$
||y_3(t)||_\a\leq Ce^{-\mu (t-t_0)}||y_3(t_0)||_\a,\Hs t\in[t_0,t_1].
$$
\el
{\bf Proof.} The proof is a slight modification of that of  \cite[Chap. II, Lemma 3.3]{Ryba}. We omit the details.  $\Box$
\Vs
%We now proceed to prove the theorem.
 Denote $G_\theta$ the local semiflow on $\~U$ generated by \eqref{rt13}. Then \eqref{rt11} implies that $W:=\~U_1\oplus \~U_2$ is a local invariant manifold of $G_\theta$ for all $\theta\in[0,1]$. Further by \eqref{rt11} and \eqref{rt13} we deduce  that  $G_\theta|_W$ is actually  independent of $\theta$.

On the other hand, if $K\subset \~U$ is an  invariant set of $G_\theta$, then by Lemma \ref{rtl1} we find that $K\subset W$. Using this fact and the independence of $G_\theta|_W$ upon $\theta$ it is easy to see that %$K\subset \~U$ is  an  invariant set of $G_{\theta_0}$ for some $\theta_0$ iff it is an  invariant set of
 all the local semiflows $G_{\theta}$ ($\theta\in[0,1]$) share the same invariant sets in $\~U$. Consequently  $K\subset \~U$ is an isolated invariant set of $G_{\theta_0}$ for some $\theta_0\in[0,1]$ iff it is an isolated invariant set of  $G_{\theta}$ for all $\theta\in[0,1]$.

Now we pick a positive number $\eta<\ve/2$ such  that
\be\label{rt25b}Q_\theta(\Xi(\eta))\subset \~U,\Hs\A \,\theta\in[0,1],\,\,\lam\in J_0.\ee
Let  $S\subset \Xi(\eta)$ be an isolated invariant set of $\Psi$. Then
 \be\label{e:S1}S_1=Q_1S=P^{12}S\ee
 is an isolated invariant set of $\~\Psi_1$ in $\~U$.  (The second equality in \eqref{e:S1} follows from  \ef{e:Q1} and  the fact that $S\subset M^{12}_\lam$.) Since $\~\Psi_1$ coincides with $G_1$, by what we have proved above  we see that $S_1$ is an isolated invariant set of $G_\theta$ for all $\theta\in[0,1]$. Further one can easily verify  that  $(G_\theta,S_\theta)$ is a static coninuation. Thus by Theorem\,\ref{t3.8} and \ef{rt10} we obtain that
\be\label{rt16}
\mb{Ind}(\Psi,S)=\mb{Ind}(\~\Psi_1,S_1)=\mb{Ind}(G_1,S_1)=\mb{Ind}(G_0,S_1),
\ee
where  $G_0$ is the local semiflow on  $\~U$ generated by the system
\be\label{rt14}\left\{\ba{lll}
\dot{x}+ B^{12} x=P^{12}h(x+\psi(x)),\\[1ex]
\dot{y}_3+ B^3 y_3=0.\ea\right.
\ee

 \Vs
{\em Step 3.} We now  calculate the index  $\mb{Ind}(G_0,S_1)$.
For this purpose, consider the following  homotopy of  system \eqref{rt14}:
\be\label{rt15}\left\{\ba{lll}
\dot{x}+ B^{12} x=P^{12}h(x+\psi(x)),\\[1ex]
\dot{y}_3+ B^3_\tau y_3=0,\ea\right.\Hs \tau\in[0,1],
\ee
where $ B^3_\tau=\tau  B^3+(1-\tau)I$. Since $\mb{Re}\,\sig( B^3)\geq2\de>0$ (see \ef{e:4B}), we have $$\mb{Re}\,\sig( B^3_\tau)\geq\min(2\de,1)>0,\Hs \A\,\tau\in[0,1].$$ Denote $\Pi_\tau$  the local semiflow generated by \eqref{rt15}. Then $\Pi_1=G_0$. Hence $S_1$ is an isolated invariant set for $\Pi_1$. Further repeating the same argument as above with minor  modifications  it is easy to deduce that
$S_1$ is an isolated invariant set of $\Pi_\tau$ for all $\tau\in[0,1]$. Now as in the preceding steps,  by Theorem \ref{t3.8} we get
\be\label{rt18}
\mb{Ind}(G_0,S_1)=\mb{Ind}(\Pi_1,S_1)=\mb{Ind}(\Pi_0,S_1).
\ee
\vs
{\em Step 4.} For $\tau=0$, \eqref{rt15} reads
\be\label{rt15b}\left\{\ba{lll}
\dot{x}+ B^{12} x=P^{12}h(x+\psi(x)),\\[1ex]
\dot{y}_3+ y_3=0.\ea\right.
\ee
An equivalent form of the system reads
\be\label{e:4.30}
\dot{w}+(P^{12} B+P^3)w=P^{12}h(x+\psi(x)),
\ee
where $w=x+y_3$. Take a number $a>0$ sufficiently large so that $aI+B:=B_a$ has bounded inverse $B_a^{-1}$. Rewrite
\ef{e:4.30} as
\be\label{e:4.31}
\dot{w}+(P^{12} B_a+P^3)w=ax+P^{12}h(x+\psi(x)):=\~h_a(x).
\ee
Since $L_\lam$ commutes with $P_\lam^{ij}$, by \ef{e:4.35} it is trivial to verify that $B=B_\lam:=T_\lam L_\lam T_\lam^{-1}$ commutes with $P^{ij}=P_{\lam_0}^{ij}$ for all $\lam\in J_0$. Hence $P^{12} B_a=B_aP^{12}$. Using this fact one can easily show that  $\~B_a:=P^{12} B_a+P^3$ is invertible with
$$
\~B_a^{-1}=P^{12}B_a^{-1}+P^3.
$$

Define an operator  $\cK$ on $\~U=(\~U_1\oplus\~U_2)\oplus \~U_3:=\~\Om\oplus \~U_3$ as
$$
\cK(w)=\~B_a^{-1}\~h_a(x),\Hs \A\,w=x+y_3\in \~\Om\oplus \~U_3.
$$
Then
$$
\cK(w)=(P^{12}B_a^{-1}+P^3)\~h_a(x)=P^{12}B_a^{-1}\~h_a(x)\in E^{12}.
$$
Hence $\cK:\~U\ra E$ is a finite dimensional operator. Therefore    by the definition of the equilibrium  index and the reduction property of the Leray-Schauder degree, we deduce that
$$
\mb{Ind}(\Pi_0,S_1)=\mb{deg}(I-\cK,N,0)=\mb{deg}\(I-\cK|_{\~\Om}, \~N,0\),
$$
where $N\subset \~U$ is a closed isolating neighborhood of $S_1$, and $\~N=N\cap\~\Om$. Thus by \ef{rt16} and \ef{rt18} we have
\be\label{e:4.32}
\mb{Ind}(\Psi,S)=\mb{Ind}(G_0,S_1)=\mb{Ind}(\Pi_0,S_1)=\mb{deg}\(I-\cK|_{\~\Om}, \~N,0\).
\ee

Denote $\Psi^{12}$ the local semiflow generated by  the first equation in \ef{rt15b}:
\be\label{rt19}
\dot{x}+ B^{12} x=P^{12}h(x+\psi(x)).
\ee
($\Psi^{12}$ is well-defined on $\Om=V_1\oplus V_2$.)
 As in Remark \ref{p:4.1} we know  that $S_1=P^{12}S$ is an isolated invariant set of
$\Psi^{12}$ with $\~N$ being an isolating neighborhood. Hence  we have
\be\label{e:4.32b}
\mb{Ind}(\Psi^{12},S_1)=\mb{deg}(I-\cF, \~N,0),
\ee
where
$$\cF(x)=(aI+B^{12})^{-1}[ax+P^{12}h(x+\psi(x))],\Hs x\in \~\Om.$$
On the other hand, it is trivial to verify that
$
(aI+B^{12})^{-1}=P^{12}B_a^{-1}|_{E^{12}}.
$
Therefore
$$
\cF(x)=P^{12}B_a^{-1}\~h_a(x)=\cK|_{\~\Om}(x),\Hs x\in \~\Om.
$$
Combing this with  \ef{e:4.32b} and \ef{e:4.32} one concludes that
 \be\label{rt17}
\mb{Ind}(\Psi,S)=\mb{Ind}(\Psi^{12},S_1).
\ee
\vs

{\em Step 5.} Now we pay some attention to the calculation of $\mb{Ind}(\Psi^{12},S_1)$. At this point, we are in a quite similar situation as we were at the beginning of calculating the index $\mb{Ind}(\Psi,S)$ with system \ef{rt2} replaced by  \eqref{rt19}.

Write $x=x_1+x_2$, where $x_i=P^ix$. Then \eqref{rt19} can be reformulated as
\be\lb{rt20}
\left\{\ba{lll}
\dot{x}_1+ B^{1} x_1=P^{1}h(x+\psi(x)),\\[1ex]
\dot{x}_2+ B^2 x_2=P^{2}h(x+\psi(x)).\ea\right.
\ee
Repeating a similar  argument as in Steps 1 and  2 leading to \ef{rt16} with some corresponding modifications, it can be shown that there exists $\b>0$ (independent of $\lam\in J_0$) such that
for any isolated invariant set $K$ of \ef{rt19} with $K\subset \mB^1_\b\oplus \mB^2_\b$,
\be\label{rt22}
\mb{Ind}(\Psi^{12},K)=\mb{Ind}(\pi,P^2K),
\ee
where $\pi$ is  the local semiflow generated by
\be\label{rt14b}\left\{\ba{lll}
\dot{x}_1+ B^{1} x_1=0,\\[1ex]
\dot{x}_2+ B^2 x_2=P^{2}h\(x_2+\zeta(x_2)\),\ea\right.
\ee
where   $\zeta(x_2)=\zeta_\lam(x_2)$ is the {\em invariant manifold } mapping  given in \eqref{rt8d}.

Let $$\cF(x_2)= B^2 x_2-P^{2}h\(x_2+\zeta(x_2)\).$$ Pick a closed  isolating neighborhood $N$ of $K$  with $N\subset \mB^1_\b\oplus \mB^2_\b$. Recalling  that $\mb{Re}\,\sig( B^1)<0$,  a simple calculation yields
\be\label{rt23}\ba{ll}
\mb{Ind}(\Psi^{12},K)&=\mb{Ind}(\pi,P^2K)\\[1ex]&=(-1)^{\mathfrak{m}_1} \mb{deg}(\cF,N\cap V_2,0)\\[1ex]&=(-1)^{\mathfrak{m}_1}\mb{Ind}(\Psi^2,P^2K),\ea
\ee
where $\Psi^2$ is the local semiflow on $V_2$ generated by the system
\be\label{rt24}
\dot{x}_2+ B^2 x_2=P^{2}h\(x_2+\zeta(x_2)\).
\ee
\vs
{\em Step 6.} We first take a positive number $\b'<\b$ such that $\Xi(\b')\subset \~U$, where $\~U$ is the neighborhood of $0$ given in Lemma \ref{rtl1}. Then choose a positive number  $\eta'<\eta$ ($\eta$ is the number in \eqref{rt25b})  small enough so that
\be\label{rt25c}Q_\theta(\Xi(\eta'))\subset \Xi(\b'),\Hs\A \,\theta\in[0,1],\,\,\lam\in J_0,\ee
where $Q_\theta$ is the mapping defined in \ef{e:Q}.

Let   $S\subset \Xi(\eta')$ be an isolated invariant set of $\Psi$ in $V$. Clearly
$$S_1:=Q_1S\subset \Xi(\b')\subset\~U.$$
Since $S_1=P^{12}S\subset E^{12}$ (see \ef{e:S1}), the inclusion  $S_1\subset \Xi(\b')$ implies   that
$$S_1\subset \mB^1_{\b'}\oplus \mB^2_{\b'}\subset \mB^1_\b\oplus \mB^2_\b.$$ Thus by \ef{rt17} and \ef{rt23} we obtain that
\be\label{rt27}
\mb{Ind}(\Psi,S)=\mb{Ind}(\Psi^{12},S_1)=\mb{Ind}(\Psi^{12},P^{12}S)=(-1)^{\mathfrak{m}_1}\mb{Ind}(\Psi^2,P^2S).
\ee
\vs
{\em Step 7.}  Let $U_\lam=T_\lam^{-1}V$. By continuity of $T_\lam^{-1}$ in $\lam$ one can easily verify that there exists a neighborhood $\cO$ of $0$ in $E$ such that
$\cO\subset U_\lam$ for all $\lam\in J_0$. Take a number $\gam>0$ sufficiently small such that $U:=\Xi(\gam)\subset \cO$; furthermore,
$$
T_\lam U\subset \Xi(\eta')\hs\mb{for all }\, \lam\in J_0.
$$

The reduction  equation \eqref{re2} on $\cM_\lam^2$ reads
\be\label{re2'}
\dot{x}_2+ L^2 x_2= P^2_\lam g(x_2+\xi(x_2)),\hs x_2=x_2(t)\in U_{\lam}^2:=T_\lam^{-1}V_2,
\ee
where $L^2=L_\lam^2$, $g=g_\lam$, and $\xi=\xi_\lam=T^{-1}_\lam \zeta_\lam T_\lam$. Write
$$\~g(x_2)= g(x_2+\xi(x_2)),\Hs x_2\in U_{\lam}^2.$$ Noticing that $T_\lam P^2_\lam=P^2 T_\lam$ (recall that $P^2=P^2_{\lam_0}$), we find that
$$\ba{ll}
T_\lam (P^2_\lam\,\~g) T_\lam^{-1}(y_2)&=T_\lam P^2_\lam g\(T_\lam^{-1}y_2+T_\lam^{-1}\zeta(y_2)\)\\[1ex]
&= P^2  (T_\lam g T_\lam^{-1}) \(y_2+\zeta(y_2)\)=P^2\~h(y_2),\Hs y_2\in V_2,\ea
$$
where $\~h(y_2)=h\(y_2+\zeta(y_2)\)$, and $h=h_\lam=T_\lam g T_\lam^{-1}$ is the function in equation \eqref{rt2}.
Since $ B^2= B^2_\lam=T_\lam L_\lam^2 T_\lam^{-1}$, by Proposition  \ref{p3.4} we deduce that  for any isolated invariant set $K$ of the local semiflow $\Psi^2=\Psi_\lam^2$ of   \ef{rt24}, one has
$$\mb{Ind}(\Psi^2_\lam,K)=\mb{Ind}(\cR^2_\lam ,T_\lam^{-1}K),
$$
where $\cR_\lam^2$ is the local semiflow of \ef{re2} on  $U_\lam^2$.

Now let $S\subset U$ be an isolated invariant set of $\Phi_\lam$. Set $K=T_\lam S$. Then by Proposition \ref{p3.4} and \eqref{rt27} we conclude that
\be\lb{rt37}\ba{ll}
\mb{Ind}(\Phi_\lam,S)&=\mb{Ind}(\Psi_\lam,K)=(-1)^{\mathfrak{m}_1}\mb{Ind}(\Psi^2_\lam,P^2K)\\[1ex]
&=(-1)^{\mathfrak{m}_1}\mb{Ind}(\cR^2_\lam ,T_\lam^{-1}P^2K)\\[1ex]
&=(-1)^{\mathfrak{m}_1}\mb{Ind}(\cR^2_\lam ,P^2_\lam T_\lam^{-1}K)=(-1)^{\mathfrak{m}_1}\mb{Ind}(\cR^2_\lam ,P^2_\lam S),\ea
\ee
which completes the proof of the validity of the second equality in \ef{ert}. (The last equality in \ef{rt37} is due to the fact that $T_\lam^{-1}P^2=P^2_\lam T_\lam^{-1}$.)

The first equality in \ef{ert} follows from Proposition \ref{p3.4} and the first equality in \ef{rt27}.
The proof of the theorem is complete. $\Box$
\section{An Index Formula for Bifurcating Invariant Sets near Equilibrium Points}
%In this section we give an equilibrium index formula on bifurcating invariant sets of system \eqref{e4.1} near equilibrium points in terms of the Euler characteristic numbers of Conley indices.

As in Section 5, we follow the same notations as in Section 4 and assume the function $f_\lam(x)$ in \ef{e4.1} satisfies all the regularity hypotheses  in Section 4. Suppose also that
$f_\lam(0)\equiv0$ for $\lam\in\R$. The main result is the following theorem.

%Let $\Phi_\lam$ be the local semiflow  of \ef{e4.1}.

%\subsection{Fixed-point index of bifurcation invariant set}
\bt\label{t5.1}In addition to  (H1)-(H3) (see Section 4), assume  that
\benu
\item[(H4)] for $\lam\in J_0$, we have
$$
\mb{\em Re}\,\sig_\lam^2<0\,\,(\mb{if }\lam<\lam_0),\hs \mb{\em Re}\,\sig_\lam^2>0\,\, (\mb{if }\lam>\lam_0).
$$
\eenu Suppose also that  $S_0=\{0\}$ is an isolated invariant set of $\Phi_{\lam_0}$.

Then there exist closed isolating neighborhood $N$ of $S_0$ with respect to $\Phi_{\lam_0}$ and  $\ve>0$ such that
  $\Phi_\lam$ has a maximal compact invariant set $K_\lam$ ($K_\lam$ may be void) in $N\setminus S_0$  for each $\lam\in[\lam_0-\ve,\,\lam_0+\ve]$, $\lam\ne\lam_0$. Furthermore,
  \be\label{eif}\ba{ll}
\mb{\em Ind}(\Phi_\lam,K_\lam)&=\left\{\ba{ll}\chi\(h(\Phi_{\lam_0},S_0)\)-(-1)^{\mathfrak{m}_1},\hs&\lam<\lam_0;\\[1ex]
\chi\(h(\Phi_{\lam_0},S_0)\)-(-1)^{\mathfrak{m}_1+\mathfrak{m}_2},\hs&\lam>\lam_0\ea\right.
\\[4ex]
&=(-1)^{\mathfrak{m}_1}\left\{\ba{ll}\chi\(h(\Phi_{\lam_0}^{2},S_0)\)-1,\hs&\lam<\lam_0;\\[1ex]
\chi\(h(\Phi_{\lam_0}^{2},S_0)\)-(-1)^{\mathfrak{m}_2},\hs&\lam>\lam_0,\ea\right.
\ea
\ee
where  $\Phi_{\lam}^{2}=\Phi_{\lam}|_{\cM_{\lam}^2}$, and  $\mathfrak{m}_i=\mb{\em dim}\,(X^i)$.% is the local semiflow associated with the reduced equation  \eqref{re2}.
\et

\br\label{r6.2} The interested reader is referred to \cite[Theorem 4.3]{LW} on when the bifurcating invariant set $K_\lam$ is nonvoid.
We also infer from the proof of \cite[Theorem 4.3]{LW} that $K_\lam$ is upper semicontinuous in $\lam$.
\er

\br For convenience, we will call the number  $\mathfrak{m}_2$ in Theorem \ref{t5.1} the crossing number at $\lam=\lam_0$ in case (H4) is fulfilled.
\er

\noindent{\bf Proof of Theorem \ref{t5.1}.} (1)
Pick a closed isolating neighborhood $N$ of $S_0$ with $N\subset U$, where $U$ is the neighborhood of $0$ given in Theorem \ref{rtt1}. Then by a very standard argument  (see e.g. \cite[Chap. I]{Ryba}) it can be shown that  there exists  $\ve>0$ such that $N$ is an isolating neighborhood of $\Phi_\lam$ for all $\lam\in[\lam_0-\ve,\,\lam_0+\ve]$.
Denote $S_\lam$ the maximal compact invariant set  of $\Phi_\lam$ in $N$. Then we infer from  Theorem \ref{t:4.7} that $S_\lam\subset \cM^2_\lam$ and  is an isolated  invariant set of $\Phi_\lam^2$ on $\cM_\lam^2$.

For convenience in statement, let us temporarily forget the conjugacy between $\Phi_\lam^2$ and the local semiflow $\cR_\lam^2$ of the reduction equation \ef{re2} and identify the two local semiflows.

Let $\lam\in[\lam_0-\ve,\lam_0)$. Then by (H4) we deduce that $\mb{Re}\,\sig(L_\lam^2)<0$, where $L_\lam^2$ is the linear operator in the reduction equation \ef{re2}. Hence $S_0=\{0\}$ is a repeller of $\Phi_\lam^2$. It follows by the standard Morse decomposition theory  of invariant sets that $S_\lam$ has a Morse decomposition $\sM=\{M_0,M_1\}$ with
$M_1=S_0$. Note that $M_0:=K_\lam$ is the maximal compact invariant set of $\Phi_\lam$ in $S_\lam$. By maximality of $S_\lam$ it can be easily seen that $K_\lam$ is also the maximal compact invariant set of $\Phi_\lam$ in $N\sm S_0$. (It may occur that $S_\lam=S_0$. In such a case we have $K_\lam=M_0=\emp$.)

If $\lam\in(\lam_0,\lam_0+\ve]$ then $S_0$ is an attractor of $\Phi_\lam^2$, and a parallel argument applies to show that $\Phi_\lam$ has a maximal compact invariant set $K_\lam$ in $N\sm S_0$.

\vs
(2) Let $K_\lam$ be the maximal compact invariant set of $\Phi_\lam$ in $N\sm S_0$. Since $0\not\in K_\lam$, we can pick a $\b=\b(\lam)>0$ with $\mB_{2\b}\subset N$ such that $\mB_{2\b}\cap K_\lam=\emp$, where $\mB_r$ denotes the ball in $E$ centered at $0$ with radius $r$. Then $N_\b:=N\sm \mB_\b$ is a closed neighborhood of $K_\lam$. Further by maximality of $K_\lam$ in $N\sm S_0$ we deduce that  $N_\b$ is an isolating neighborhood of $K_\lam$.

We infer from the Morse decomposition theory that $S_\lam$ is the union of $S_0$, $K_\lam$ and the connecting orbits between $S_0$ and $K_\lam$. Thus one concludes that
$\sE(S_\lam)\subset K_\lam\cup S_0.$ (Recall that $\sE(S)$ denotes the set of equilibrium points in $S$.)  Therefore
\be\label{e:6.1}
\mb{Ind}\,(\Phi_\lam,S_\lam)=\mb{Ind}\,(\Phi_\lam,K_\lam)+\mb{Ind}\,(\Phi_\lam,S_0).
\ee
%Hence \be\label{e5.1}\mb{Ind}\,(\Phi_\lam,K_\lam)=\mb{Ind}\,(\Phi_\lam,S_\lam)-\mb{Ind}\,(\Phi_\lam,S_0).\ee

On the other hand, by homotopy property of the equilibrium index we have
$$
\mb{Ind}\,(\Phi_\lam,S_\lam)=\mb{Ind}\,(\Phi_{\lam_0},S_{\lam_0})=\mb{Ind}\,(\Phi_{\lam_0},S_{0}).
$$
Hence by \ef{ert} and  Proposition \ref{p3.4} we deduce that
\be\label{e:6.2}\ba{ll}
\mb{Ind}\,(\Phi_\lam,S_\lam)&=\mb{Ind}\,(\cR_{\lam_0}^{12},P_{\lam_0}^{12}S_0)=\mb{Ind}\,(\cR_{\lam_0}^{12},S_0)\\[1ex]
&=\chi\(h(\cR_{\lam_0}^{12},S_0)\)=\chi\(h(\Phi_{\lam_0}^{12},S_0)\)\\[1ex]
&=(\mb{by \ef{erth}})\\[1ex]
&=\chi\(\Sig^{\mathfrak{m}_1}\wedge h(\Phi_{\lam_0}^2,S_0)\)=\chi\(h(\Phi_{\lam_0},S_0)\).\ea
\ee
It then follows by \ef{e:6.1}  that
\be\label{e:6.3}
\mb{Ind}\,(\Phi_\lam,K_\lam)=\chi\(h(\Phi_{\lam_0},S_0)\)-\mb{Ind}\,(\Phi_\lam,S_0).
\ee

Simple calculations  show that
\be\label{e5.2}
\mb{Ind}\,(\Phi_\lam,S_0)=\left\{\ba{ll}(-1)^{\mathfrak{m}_1},\hs&\lam<0;\\[1ex]
(-1)^{\mathfrak{m}_1+\mathfrak{m}_2},\hs&\lam>0.\ea\right.
\ee
Therefore
$$
\mb{Ind}\,(\Phi_\lam,K_\lam)=\left\{\ba{ll}\chi\(h(\Phi_{\lam_0},S_0)\)-(-1)^{\mathfrak{m}_1},\hs&\lam<0;\\[1ex]
\chi\(h(\Phi_{\lam_0},S_0)\)-(-1)^{\mathfrak{m}_1+\mathfrak{m}_2},\hs&\lam>0.\ea\right.
$$

We also infer from \ef{e:6.2} that
$$\ba{ll}
\chi\(h(\Phi_{\lam_0},S_0)\)=\chi\(\Sig^{\mathfrak{m}_1}\wedge h(\Phi_{\lam_0}^{2},S_0)\)
=(-1)^{\mathfrak{m}_1}\chi\(h(\Phi_{\lam_0}^{2},S_0)\).\ea
$$
Combining this with \ef{e:6.1} and \ef{e5.2} one immediately concludes the validity of the second equality in \eqref{eif}. $\Box$
%$$\mb{Ind}\,(\Phi_\lam,K_\lam)=\left\{\ba{ll}(-1)^{\mathfrak{m}_1}\[\chi\(h(\Phi_{\lam_0}^{2},S_0)\)-1\],\hs&\lam<0;\\[1ex]
%(-1)^{\mathfrak{m}_1}\[\chi\(h(\Phi_{\lam_0}^{2},S_0)\)-(-1)^{\mathfrak{m}_2}\],\hs&\lam>0.\ea\right.$$
%This  completes the proof of the theorem. $\Box$

%By the  homotopy invariance of the Conley index (see \cite{Ryba}), we have
%$$h(\Phi_\lam,S_\lam)=const.$$
%for $\lam\in[-\ve,\ve]$. Now if $h(\Phi_0,S_0)\ne \Sig^{\mathfrak{m}_1}$,  then since $h(\Phi_\lam,S_0)=\Sig^{\mathfrak{m}_1}$ for $\lam\in[-\ve,0)$, we deduce that $S_\lam\ne S_0$.   We infer from the reduced equation \ef{re2} that $S_0$ is a repeller of $\pi=\Phi_\lam|_{\cM_\lam^2}$. Hence $S_\lam$ has an attractor-repeller pair $(K_\lam,S_0)$. It is easy to see that $K_\lam$ is the  maximal compact invariant set in $N\sm S_0$. Thus (1) holds true.

\section{A Global Static  Bifurcation Theorem}
%In this section we first prove an existence result for local bifurcation branch.

%Consider the equation \be\label{sp1} Au=f_\lam(u),\Hs u\in E=X^\a,\ee or equivalently
%\be\label{sp2} L_\lam u=g_\lam(u),\Hs u\in E, \ee where $L_\lam$ and $g_\lam$ is a s in Section 4.
We follow the same notations in the preceding sections. Furthermore, we assume the function $f_\lam(x)$ in \ef{e4.1} satisfies all the regularity hypotheses  in Section 4.

Let $\sX=E\X\R$. $\sX$ is equipped with metric $d$ defined as
$$
d\((u,\lam),\,(v,\lam')\)=||u-v||_\a+|\lam-\lam'|, \Hs \A\, (u,\lam),\, (v,\lam')\in\sX.
$$
Given $\cZ\subset \sX$ and $\lam\in \R$, denote $\cZ[\lam]$ the $\lam$-section of $\cZ$:
$$
\cZ[\lam]=\{u:\,\,(u,\lam)\in \cZ\}.
$$

Let $\Phi_\lam$ be the local semiflow of \eqref{e4.1}, and denote $\Pi$  the {\em skew-product flow} of the family $\Phi_\lam$ ($\lam\in \R$) on $\sX$,
$$
  \Pi(t)(x,\lambda)=(\Phi_\lambda(t)x,\lambda), \Hs (x,\lambda)\in\sX,\,\,t\geq0.
  $$
%\be\label{e6.1}
%\Pi(t)(u,\lam)=\(\Phi_\lam(t)u,\,\lam\),\Hs\A\,(u,\lam)\in\sX.
%\ee

As in the  proof of Theorem \ref{t3.5} it can be shown that  the set
$\{f_\lam(u):\,\,(u,\lam)\in\cB\}$ is bounded in $X$ for any bounded subset $\cB$ of $\sX$. Using this simple  fact and applying some  fundamental theory on abstract evolution equations (see e.g.\,\,\cite[Chap. 3]{Henry} and \cite[Chap. I, Theorem 4.4]{Ryba}), it can be shown by very standard argument that $\Pi$ is {\em asymptotically compact}.
Consequently   each bounded closed invariant set $\cK$ of $\Pi$ is necessarily compact.

Assume  $f_\lam(0)\equiv0$ for all $\lam\in\R$, hence  $u=0$ is always a trivial equilibrium of $\Phi_\lam$. Given $\cU\subset \sX$, denote
$$\ba{ll}
\sC(\cU)=\mb{cl}_\sX\,\{(u,\lam)\subset \cU:\,\, u \mb{ is a nontrivial equilibrium of }\Phi_\lam\}.\ea
$$
%where ``nontrivial\,'' means $u\ne0$.

\begin{definition} Let $\cU\subset \sX$ be a closed neighborhood of  $(0,\mu)$.
The (static) {bifurcation branch} $\Gamma_\cU(0,\mu)$ of $\Phi_\lam$ from  $(0,\mu)$ in $\cU$
 %denoted by $\Gamma_\cU(0,\lam_0)$,
 is defined to be the   component of ${\sC(\cU)}$ which contains $(0,\mu)$.
%; see Figure 5.1.
\end{definition}

\subsection{Existence of local bifurcation branch}
Let us first give an existence result of a nontrivial local bifurcation branch.

Let $\cM_\lam^2$ and $\cM_\lam^{12}$ be the local invariant manifolds of $\Phi_\lam$ given in Theorem \ref{t:4.7}, and denote $\Phi_\lam^2=\Phi_\lam|_{\cM_\lam^2}$, $\Phi_\lam^{12}=\Phi_\lam|_{\cM_\lam^{12}}$. 
\bt\label{t6.2}Assume  (H1)-(H4) (see Section 4 and Theorem \ref{t5.1})  are fulfilled. Suppose $S_0=\{0\}$ is an isolated invariant set of $\Phi_{\lam_0}$, where $\lam_0$ is the number appearing in (H4).
Let $N$ be the isolating neighborhood of $S_0$ given in Theorem \ref{t5.1}.
Then there exists $\ve>0$ such that the following assertions hold.
\benu
\item[$(1)$]If $\chi\(h(\Phi_{\lam_0}^{2},S_0)\)\ne (-1)^{\mathfrak{m}_2}$ then
  $ \Gamma[\lam_0+\ve]\ne\emp$.
\item[$(2)$]  If
$\chi\(h(\Phi_{\lam_0}^{2},S_0)\)\ne 1$ then
$ \Gamma[\lam_0-\ve]\ne\emp.$
\eenu
Here $\Gamma=\Gamma_\cN(0,\lam_0)$, and $\cN=N\X[\lam_0-\ve,\,\lam_0+\ve]$.
\et

\Vs\noindent{\bf Proof.} For simplicity, we set $\lam_0=0$.
Choose an  $\ve>0$ such  that the assertions in Theorem \ref{t5.1} %  (\ref{3.28}) and (\ref{3.29})
hold. Let $S_\lam$ be the maximal compact invariant set of $\Phi_\lam$ in $N$, and $K_\lam$ the maximal compact invariant set of $\Phi_\lam$ in $N\sm S_0$.  We may restrict $\ve$ small enough so that $N$ is an isolating neighborhood of $S_\lam$ for all $\lam\in[-\ve,\ve]$.

For definiteness we assume $\chi\(h(\Phi_0^{2},S_0)\)\ne (-1)^{\mathfrak{m}_2}$ and prove that
\be\label{sp4}
\Gamma[\ve]\ne\emp.
\ee

Let us  first show that for any $0<\mu<\ve$, there is  a component $\cZ$  of $\sC(\cN_\mu)$, where $\cN_\mu=N\X[\mu,\ve]$, such that
\be\label{e:4.10}\ba{ll}
\cZ[\mu]\ne\emp\ne \cZ[\ve].\ea
\ee
Set $\Lam=[\mu,\ve]$, and let $\cK=\Cup_{\lam\in\Lam}K_\lam\X\{\lam\}$. Then $\cK$ is a bounded invariant set of the skew-product flow $\Pi$. By upper semicontinuity of $K_\lam$ in $\lam$ (see Remark \ref{r6.2}) one can easily verify that $\cK$ is closed. Remark \ref{r2.1} then asserts  that $\cK$ is compact in $\sX$. As $(0,\lam)\not\in \cK$ for  $\lam\in\Lam$,  we have
\be\label{lbc}
\min_{\lam\in\Lam}d(0,K_\lam):=r_0>0.
\ee
It is also trivial to deduce that $\min_{\lam\in\Lam}d(K_\lam,\pa N):=r_N>0$.
 Pick a  number  $0<r<r_0$. Then $\Om:=N\sm \mB_r$ is a closed  isolating neighborhood of $K_\lam$  for all $\lam\in\Lam$.

 Take  a number $a>0$ sufficiently large so  that $\mb{Re}\,\sig(aI+A)>0$. Set $$F_\lam=(aI+A)^{-1}(aI+f_\lam).$$ As in the proof of Theorem \ref{t3.5} we know that $\Cup_{\lam\in \Lam}F_\lam(B)$ is precompact in $E$ for any bounded subset $B$ of $E$.

 We infer from Theorem \ref{t5.1} that
 $$
\deg(I-F_\mu,\Om,0)=\mb{Ind}(\Phi_\mu,K_\mu)=\chi\(h(\Phi_{0}^{2},S_0)\)-(-1)^{\mathfrak{m}_2}\ne 0.
$$
This allows us to apply  the classical Leray-Schauder continuation theorem (see e.g. Mawhin \cite[Section 2]{Mawh1}) to deduce that there is a component $\cZ\in\sC(\Om\X\Lam)$ such that $\cZ[\mu]\ne\emp\ne \cZ[\ve]$, from which \ef{e:4.10} immediately follows.

We are now ready to complete the proof of the theorem. Take a sequence of positive numbers  $\mu_k\ra 0$. For each $\mu_k$,  pick  a connected component $\cZ_k$ of $\sC(\cO_{\mu_k})$ such that
$
\ba{ll}\cZ_k[\mu_k]\ne\emp\ne \cZ_k[\ve].\ea
$
By Lemma \ref{l:2.3} we may assume  that $\cZ_k$ converges in the sense of Hausdorff distance $\de_{\mbox{\tiny H}}(\.,\.)$ to a compact set $\cZ_0$. Then  $\cZ_0$ is a continuum in $\sC(\cN)$ with $\cZ_0[0]\ne\emp\ne\cZ_0[\ve]$. On the other hand, by the choice of $N$ we have  $\sC(\cN)[0]=\{0\}$. Therefore one concludes that  $0\in \cZ_0$.\, $\Box$

\subsection{Global static bifurcation theorem}
For simplicity, in this subsection we rewrite $\mB_{E}(0,r)=\mB_r$\,. Given $\cC\subset \sX$ and $J\subset \R$, we also denote
$$
\cC|_J=\{(u,\lam)\in\cC:\,\,\lam\in J\}.
$$

Let $M_0=(0,\lam_0)$ be a bifurcation point.

\begin{definition}The global (static) bifurcation branch $\Gamma$ of $M_0$ is defined to be the   bifurcation branch of $M_0$ in $\sX$.

The right-hand side  global  bifurcation branch $\Gamma^+$ of $M_0$ is defined as
$$\ba{ll}
\Gamma^+=\lim_{\ve\ra0}\Gamma^+_{\ve}=\Cup_{\ve>0}\Gamma^+_{\ve},\ea
$$
where $\Gamma^+_\ve$ denotes the   bifurcation branch of $M_0$ in $\sX\sm \(\mb{\em B}_\ve\X(-\ve,0)\)$.

Similarly one can define the  left-hand side global  bifurcation branch $\Gamma^-$.

\end{definition}

\br
One  easily verifies that both $\Gamma^-$ and $\Gamma^+$ are connected; moreover, we have  $\Gamma=\Gamma^+\cup\Gamma^-$.
\er

Our main result in this section is  the following global bifurcation theorem.
\bt\label{gbt} Assume  (H1)-(H4)   are fulfilled. Suppose  $S_0=\{0\}$ is an isolated invariant set of $\Phi_{\lam_0}$; furthermore,
 $$\chi\(h(\Phi_{\lam_0}^{2},S_0)\)\ne (-1)^{\mathfrak{m}_2}\mb{ or } 1.$$
Then one of the following cases  occurs (see Figures 6.3-6.5).
\begin{enumerate}
\item[$(1)$]  $\Gamma$ is unbounded.
%\item[(1)] There exist closed isolating neighborhood $N$ of $0$ and $\ve>0$ with $\cN:=N\X[-\ve,\ve]\subset \Om$ such that $\Gamma_\cN:=\Gamma_\cNM_0\subset \Gamma$; moreover,
% $$ \Gamma_\cN\cap (N\X\{-\ve\})\ne\emp\ne \Gamma_\cN\cap (N\X\{\ve\}).$$
\item[$(2)$] There exists $\lam_1\ne\lam_0$  such that $(0,\lam_1)\in\Gamma$\,.
\item[$(3)$]  $ \Gamma^+\cap \Gamma^-\ne \{M_0\}$, in which case  both  $\Gamma^-$ and $\Gamma^+$  return back to  $M_0$.\end{enumerate} \et

\begin{center}\includegraphics[width=4cm]{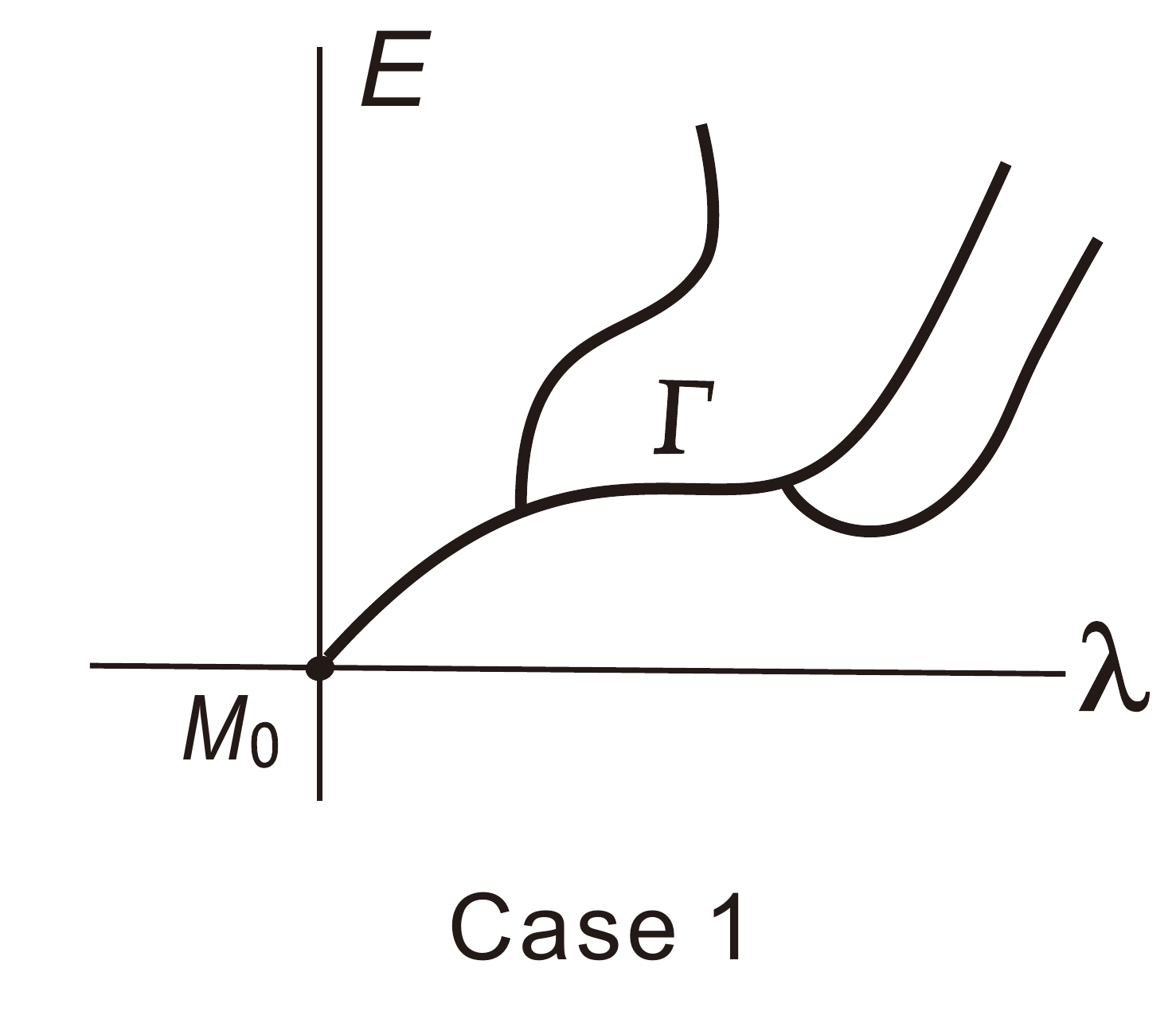} \hspace{2cm} \includegraphics[width=4cm]{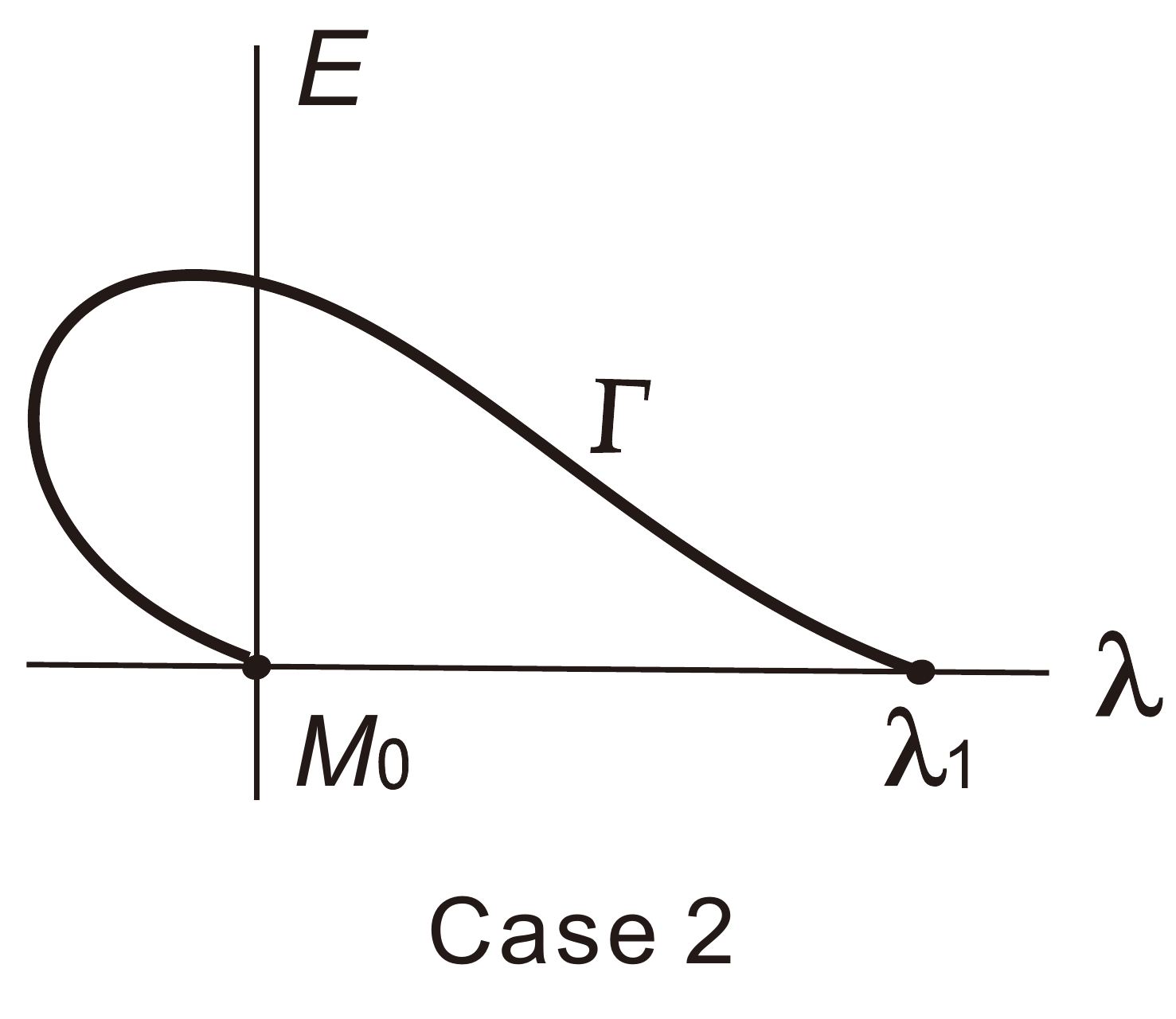}\end{center}
\begin{center}Figure 6.3: \,$\Gamma$ is unbounded. \hspace{1cm} Figure 6.4: \, $\Gamma$ connects to $(0,\lam_1)$. \end{center}
\vs
\begin{center}\includegraphics[width=4cm]{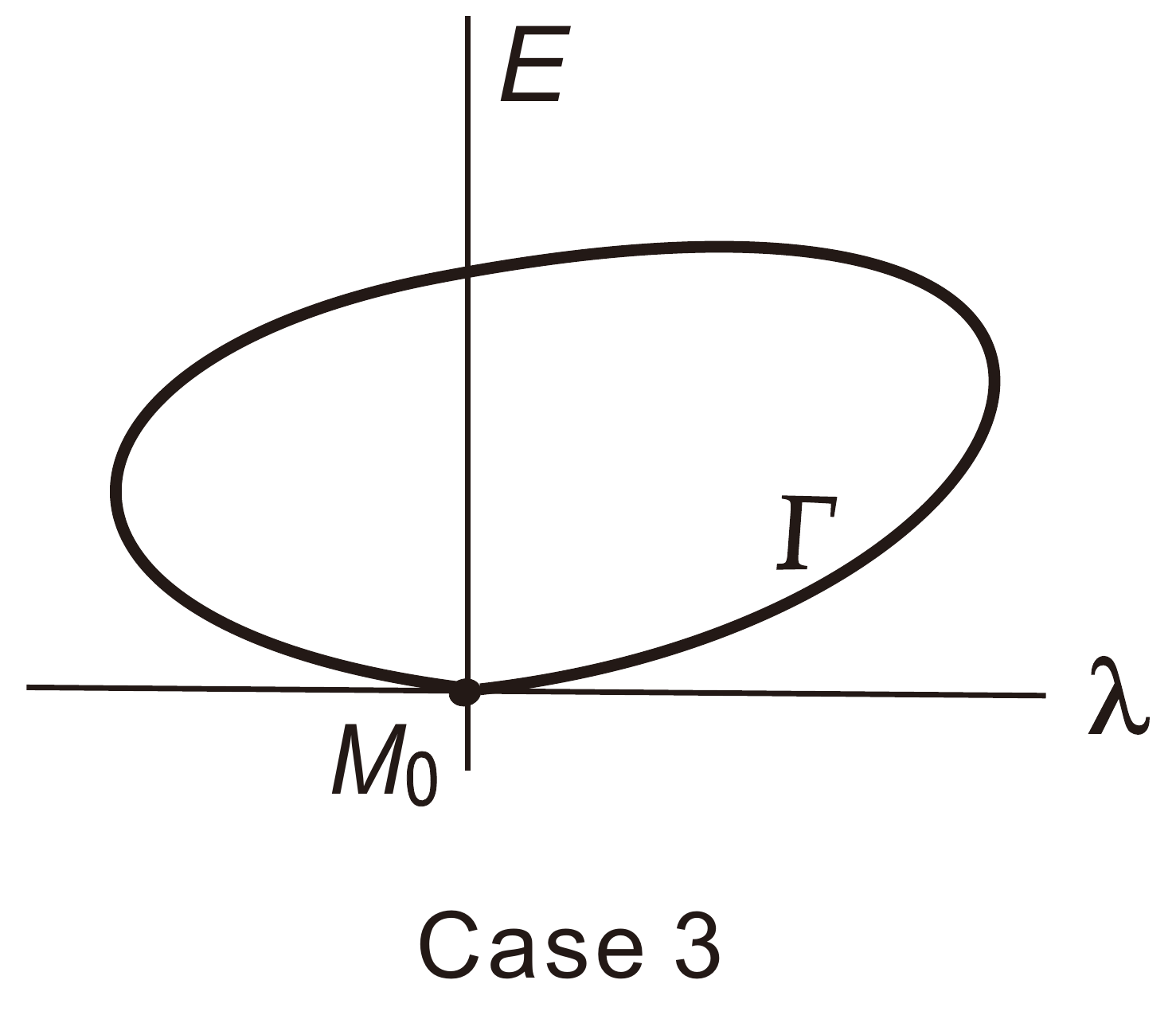} %\hspace{1.2cm} \includegraphics[width=5cm]{fig6-6.pdf}
\end{center}
\centerline { Figure 6.5: \,$\Gamma^\pm$ return back to $M_0$. }
\Vs
\noindent{\bf Proof.} We may assume $\lam_0=0$.
We argue by contradiction and suppose that none of the cases  (1)-(3) occurs. Then both $\Gamma^\pm$ are
bounded closed subsets of $\sX=E\X\R$. Moreover,
\be\label{gba1}
\Gamma^\pm\cap (\{0\}\X \R)=\{M_0\},
\ee
\be\label{gba2}
 \Gamma^+\cap\Gamma^-=\{M_0\};
\ee
see Fig. \ref{fig7-4}.  By asymptotic compactness of  $\Pi$ we know that  $\Gamma^\pm$ are compact.

\vs
{\em Step 1.} %By the assumption that $\Gamma^+\cap \Gamma^-=\{M_0\}$, one can easily verify that the point $M_0$ is separated in both $\Gamma^+|_{\R^-}$ and $\Gamma^-|_{\R^+}$.
By  Theorem \ref{t6.2} we have either $\Gamma^+\sm\{M_0\}\ne\emp$, or $\Gamma^-\sm\{M_0\}\ne\emp$. Let $$\cR=\Gamma^+|_{\R^-}\sm\{M_0\},\hs   \cL=\Gamma^-|_{\R^+}\sm\{M_0\}.$$
We may assume $\cR\ne\emp\ne  \cL.$ In the case when  $\cR=\emp$ or $ \cL=\emp$, the argument  is simpler and can be obtained by slightly modifying  the one below.
In the following we first show that $\cR$ and $ \cL$ are closed and hence are compact. For this purpose it suffices to check that  $M_0$ is an isolated point in both the sets  $\Gamma^+|_{\R^-}$ and $\Gamma^-|_{\R^+}$.

We argue by contradiction and suppose that $M_0$ is not  isolated in, say,   $\Gamma^+|_{\R^-}$. Then
there exists a sequence
%\be\label{sp13}d\(0,\, \cR[0]\)>0,\hs d\(0, \cL[0]\)>0.\ee
   $(u_n,\lam_n)\in \cR$ converges to $M_0$. Since $u=0$ is an isolated equilibrium of $\Phi_0$, it is clear that
   $\lam_n\ne 0$ (hence $\lam_n<0$) for all $n$ sufficiently large.

   Pick a closed statically isolating neighborhood $V$ of $0$ with respect to $\Phi_0$. There is  $\ve>0$ such that $V$ is statically isolating with respect to $\Phi_\lam$ for  $\lam\in[-\ve,\ve]$.
 It can be assumed  that $\lam_n<0$ and  $(u_n,\lam_n)\in \cV:=V\X[-\ve,\ve]$  for all $n$.

Denote $\cC_n$ the continuum of $\cV\cap\Gamma^+|_{\R^-}$ containing $(u_n,\lam_n)$. If $M_0\in \cC_n$ then by definition we find that $\cC_n\subset \Gamma^-$. Hence $(u_n,\lam_n)\in \Gamma^+\cap\Gamma^-$, which contradicts \ef{gba2}. Thus we deduce that $M_0\not\in \cC_n$. As $u=0$ is the only equilibrium of $\Phi_0$ in $V$,  it follows that  $\cC_n\cap (V\X \{0\})=\emp.$ Further, because $\Gamma^+\cap (\pa V\X[-\ve,\ve])=\emp$, by connectedness of $\Gamma^+$ one concludes that   $\cC_n[-\ve]\ne\emp$.

By Lemma \ref{l:2.3} we may assume  that $\cC_n$ converges in the sense of Hausdorff distance  to a compact subset $\cC_0$. Then $\cC_0$ is a continuum in $\Gamma^+$. Clearly $M_0\in \cC_0$ and $\cC_0[-\ve]\ne\emp$. Hence $\cC_0[\lam]\ne\emp$ for all $\lam\in[-\ve,0]$. We infer from  \ef{gba1}   that
\be\label{ca}0\not\in \cC_0[\lam],\Hs \lam\in[-\ve,0).\ee
 Thereby $\cC_0\subset\Gamma^-$. But this and \ef{ca} contradict \ef{gba2}.

\vs
{\em Step 2.}  By \ef{gba1} and the definition of $\cR$ and $\cL$ it is clear  that  $(0,\lam)\not\in \cR\cup \cL$ for all $\lam\in\R$. Thus by compactness of $\cR$ and $\cL$, there exists $\de_0>0$ such that
\be\label{sp3a}||u||_\a\geq 3\de_0>0,\Hs\A\,(u,\lam)\in\cR\cup\cL.
\ee

Let $N\subset E$ and $\ve>0$ be given  as  in  Theorem \ref{t5.1}. We may assume $N$ is chosen sufficiently small so that $N\subset \mB_{\de_0}$.
Clearly \be\label{spa1}\ba{ll}\mB(\cR[0],\de_0)\cap  N=\emp=\mB(\cL[0],\de_0)\cap  N;\ea\ee
 see Fig. \ref{}.
Pick a number
$0<\de_1<\de_0$ with $\mB_{4\de_1}\subset N$. Denote $S_\lam$  the maximal invariant set of $\Phi_\lam$ in $N$.
Since \be\label{e5.11}\lim_{\lam\ra0}d_{\mbox{\tiny H}}(S_\lam,S_0)=0,\ee there exists  $0<a\leq\ve$ such that
\be\label{sp3c}
\sE_{\Phi_\lam}(N)\subset S_\lam\subset \mB_{\de_1}\subset \mB_{4\de_1}\subset N,\Hs \lam\in [-a,a]:=J_a\,.
\ee

By  compactness of $\Gamma^\pm$ it is easy to verify  that $\Gamma^\pm[\lam]$ are upper semicontinuous in $\lam$.  Thus  we can restrict    $a$ sufficiently small so that
\be\label{sp8}\ba{ll}
\Gamma^+[\lam]\subset \mB(\Gamma^+[0],\de_1)=\mB_{\de_1}\cup \mB(\cR[0],\de_1),\Hs \lam\in J_a,\ea
\ee
\be\label{sp8b}\ba{ll}
\Gamma^-[\lam]\subset \mB(\Gamma^-[0],\de_1)=\mB_{\de_1}\cup \mB( \cL[0],\de_1),\Hs \lam\in J_a.\ea
\ee
(Note  that $\Gamma^+[0]=S_0\cup \cR[0]$, and $\Gamma^-[0]=S_0\cup  \cL[0]$.) 
%Let $\cN:=N\X J_a$, and set
 %$$\cA=\Gamma^+|_{(-\8,\,a]} \sm \cN,\hs  \cB=\Gamma^-|_{[-a,\,\8)}\sm \cN.$$ Then $\cR$ and $ B$ are compact.
Let 
$$\ba{ll}
\cR_a=\cR\cup\(\Gamma^+\cap \(\,\ol\mB(\cR[0],\de_1)\X J_a\)\),\hs \cL_a=\cL\cup\(\Gamma^-\cap \(\,\ol\mB(\cL[0],\de_1)\X J_a\)\).\ea$$
By  \ef{sp3a}, \ef{sp8}, \ef{sp8b} and the choice of $\de_1$ we clearly have
\be\label{sp10a}
||u||_\a>2\de_0>2\de_1,\Hs\A\,(u,\lam)\in\cR_a\cup\cL_a.\ee
%\be\lb{sp10b}\,\,\,\,||u||_\a>2\de_0>2\de_1,\Hs\A\,u\in B[\lam],\,\,\lam\geq -a.\ee
Thus by  \ef{gba2} one finds  that
 $\cR_a\cap \Gamma^-=\emp$ and $ \cL_a\cap \Gamma^+=\emp.$ The choice of $\de_1$ and \ef{sp3a} also imply 
 $$
 d(\cR_a,\cN)>0,\hs d(\cL_a,\cN)>0,
 $$ 
 where $\cN:=N\X J_a$.
  Therefore we can pick a  $0<{\de_2}<\de_1/2$ such that
\be\label{sp11a}
 d(\cR_a,\Gamma^-)>4{\de_2},\hs d( \cL_a, \Gamma^+)>4{\de_2},
\ee
and
 \be\label{sp11b}
 d(\cR_a,\cN)>4{\de_2},\hs d(\cL_a, \cN)>4{\de_2}.
\ee

It is easy to check that 
\be\label{sp11c}
\cR_a=\Gamma^+|_{(-\8,a]}\sm \cN,\hs \cL_a=\Gamma^-|_{[-a,\8)}\sm \cN.
\ee
Let  $\cZ:=\Gamma\cap \cN$. Then since $\Gamma=\Gamma^+\cup\Gamma^-$, it follows from \ef{sp11c} that 
\be\label{sp20}
\Gamma|_{J_a}=\cZ\cup \cR_a|_{J_a}\cup \cL_a|_{J_a}.
\ee

\vs
{\em Step 3.}
By \ef{e5.11} there is  $c>0$ with $c<\frac{1}{2}\min\{\de_2,a\}$  such  that
\be\label{sp12}
\sE_{\Phi_\lam}(N)\subset S_\lam\subset \mB_{\de_2},\Hs \A\,\lam\in J_c:=[-c,c],
\ee
where $S_\lam$ is  the maximal invariant set of $\Phi_\lam$ in $N$.

Let $A_{c}:=\Cup_{|\lam|\geq c/2}\Gamma[\lam]$. Then by compactness of $\Gamma$ we deduce that  $A_{c}$ is a compact subset of $E$. As $0\not\in A_{c}$, one has
\be\label{e6.8}d(0,A_{c}):=4\kappa>0.\ee
Take a number   $0<r<\frac{1}{2}\min\{c,\kappa\}.$
%Utilizing the Separation Lemma,
By some standard argument (see e.g. the proof of \cite[Theorem 6.2]{LW}, one can find a bounded closed
neighborhood $\cO$ of $\Gamma$ with $\cO\subset \mB_\sX(\Gamma,r)$ %, where  $\mB_\sX(\Gamma,r)$ is the $r$-neighborhood of $\Gamma$ in $\sX$,
such that
\be\label{5.10}\ba{ll}
\sC(\sX)\cap \pa\cO=\emp.\ea
\ee

%Set $$\cF=\cO\cap\cN,\hs \cG=\cO|_{(-\8,a]}\sm \cN,\hs \cH=\cO|_{[-a,\8)}\sm \cN.$$
Let $$\sX_0=E\X J_a,\hs \sX_1=E\X(-\8,a], \hs\mb{and }\,\sX_2=E\X[-a,\8).$$ Define
\be\label{spfn}
\cF=\cO\cap \ol\mB_{\sX_0}(\cZ,\de_2),\hs
\cG=\cO\cap \ol\mB_{\sX_1}(\cR_a, \de_2),\hs \cH=\cO\cap \ol\mB_{\sX_2}( \cL_a, \de_2).
\ee
%where  $\mB_{\sX_i}(\cK, \eta)$  denotes the $\eta$-neighborhood of $\cK$ in $\sX_i$.
We infer from \ef{sp11a} and   \ef{sp11b}  that $\cF$, $\cG$ and $\cH$ are disjointed.  
 Recalling that $\cO\subset \mB_{\sX}(\Gamma,r)$ and $r<c/2<\de_2/4$, one trivially verifies that 
 \be\label{cgh}
\cG\subset  \mB_{\sX_1}(\cR_a, \de_2),\hs \cH\subset \mB_{\sX_2}( \cL_a, \de_2).\ee
Furthermore, the following basic facts hold true.
\bl\label{l7.2} Let $
F_\lam=\cF[\lam],$ $G_\lam=\cG[\lam],$ and  $H_\lam=\cH[\lam]$. Then 
%\item[$(2)$]
\be\label{sp26}
F_\lam\subset \ol{\mb{\em\mB}}_{2\de_1}\subset \mb{\em\mB}_{4\de_1}\subset N, \Hs \A\,\lam\in J_a,
\ee
%\item[$(3)$] \,
\be\label{sp21}
\cO[\lam]=F_\lam\cup G_\lam\cup H_\lam,\Hs \lam\in J_c,
\ee
\be\label{sp9g}
\ba{ll}||u||_\a\geq 3\de_2 ,\Hs \A\,(u,\lam)\in \cG\cup\cH, %\(\Cup_{\lam\leq a}G_\lam\)\Cup \(\Cup_{\lam\geq-a}H_\lam\),
\ea
\ee
where  $\mb{\em\mB}_r$ denotes the ball in $E$ of center $0$ and radius $r$.
%\eenu
\el
{\bf Proof.} (1)\, As $\de_2<\de_1/2$, the inequality  \ef{sp9g} directly follows from  \ef{sp10a}.

\vs
(2) \,Note that  $\cZ[\lam]\subset S_\lam$ for all $\lam\in J_a$. Thus  by \ef{sp3c} we have $\cZ\subset \mB_{\de_1}\X J_a$. Observing that  if  $||u||_\a> 2\de_1$ then 
$$
d((u,\lam),\cZ)\geq d\((u,\lam),\mB_{\de_1}\X J_a\)=d(u,\mB_{\de_1})> \de_1>\de_2,\Hs \lam\in J_a,
$$
one  concludes  that
$
\ol\mB_{\sX_0}(\cZ,\de_2)\subset \ol\mB_{2\de_1}\X J_a,
$
from which \eqref{sp26} immediately follows.

%Now we verify the validity of \ef{sp23}.Suppose the contrary. Then one would find a sequence $(u_k,\lam_k)$ with $u_k\in \pa F_{\lam_k}$  such that  $(u_k,\lam_k)\ra M_0$. Noticing that$(u_k,\lam_k)\in\pa\cO$, one finds  that $M_0\in \pa\cO$, which leads to a contradiction!

\vs (3) It is obvious that $\cO[\lam]\supset F_\lam\cup G_\lam\cup H_\lam$. So we only need to verify the converse inclusion   for each $\lam\in J_c$. For this purpose, it suffices to check that
if  $(u,\lam)\in\cO$ and $|\lam|\leq c$ then $(u,\lam)\in \cF\cup\cG\cup\cH$.

Indeed, let  $M=(u,\lam)\in\cO$, $|\lam|\leq c$. Then for any $M'=(u',\lam')\in\sX$ with $|\lam'|>a$, we have $$d(M,M')\geq |\lam'-\lam|\geq a-|\lam|\geq a-c>c>r.$$ Therefore since $d(M,\Gamma)<r$, we deduce that 
$
d(M,\Gamma)=d\(M,\Gamma|_{J_a}\).
$
Thus by \ef{sp20} one has 
$$\ba{ll}
&\min\{d\(M,\cZ\),\,d\(M, \cR_a|_{J_a}\),\,d\(M,  \cL_a|_{J_a}\)\}\\[1ex]
=&d\(M,\Gamma|_{J_a}\)=d(M,\Gamma)<r<c/2<\de_2/4<\de_2,\ea
$$
which implies  $$M\in \mB_{\sX_0}(\cZ,\de_2)\cup \mB_{\sX_1}(\cR_a, \de_2)\cup \mB_{\sX_2}( \cL_a, \de_2).$$
Hence by  \ef{spfn} we see that  $M\in \cF\cup\cG\cup\cH$. This  finishes the proof of \ef{sp21}.   $\Box$

\Vs
%We now proceed to prove Theorem \ref{gbt}.

\Vs
{\em Step 4.}
For the sake of definiteness, we assume
\be\label{sp7}\chi\(h(\Phi_0^{2},S_0)\)\ne (-1)^{\mathfrak{m}_2}\ee
and focus our attention on the interval $(0,\8)$.
Because $\cF$ is a neighborhood of $M_0:=(0,\lam_0)=(0,0)$ in $\sX$, there exist
 $\sig,\eta>0$  such that $\mB_{\sig}\X J_\eta\subset \cF,$ i.e.,  \be\label{sp23}\mB_{\sig}\subset F_\lam,\Hs \A\,\lam\in J_{\eta}:=[-\eta,\eta].\ee
By (\ref{e5.11}) we can find  a  number
$0<\mu\leq\frac{1}{2}\min\{\eta,c\}$ such that
\be\label{e:4.15} K_\lam\subset S_\lam\subset
\mB_\sig\subset F_\lam\,,\Hs\,\A \lam\in(0, 2\mu],\ee where $K_\lam$ is the bifurcating  invariant set given in Theorem \ref{t5.1}.

Let $K=\Cup_{\mu\leq\lam\leq c}K_\lam$. Then as in \eqref{lbc} we deduce that
$$
d(0,K):=4\kappa'>0.
$$

Take a number $\Lam>0$  such that 
\be\label{eol}\cO\subset E\X [-\Lam+1,\Lam-1].\ee Let $0<\rho<\min(\kappa',\kappa)$, where $\kappa$ is the number given in \eqref{e6.8}. Let  $\sY=E\X[\mu,\Lam]$, and define
$$\ba{ll}
\cV=\cO|_{[\mu,\Lam]}=\cO\cap \sY,\hs
\cW=\cV\sm\(\mB_\rho\X[\mu,\Lam]\).\ea
$$ Clearly $\cV$ is closed in $\sY$. Since $\mB_\rho\X[\mu,\Lam]$ is open in $\sY$, we see that $\cW$ is closed in $\sY$ as well.
We claim that \be\label{5.16}\sC(\cW)=\sC(\cV):=\sC.\ee To see this, by definition it suffices to show that if $\lam\in [\mu,\Lam]$, then  any equilibrium $e\ne0$ of $\Phi_\lam$ in
$\cV[\lam]$ is  contained in $\cW[\lam]$.

\begin{center}\Hs\includegraphics[width=6.6cm]{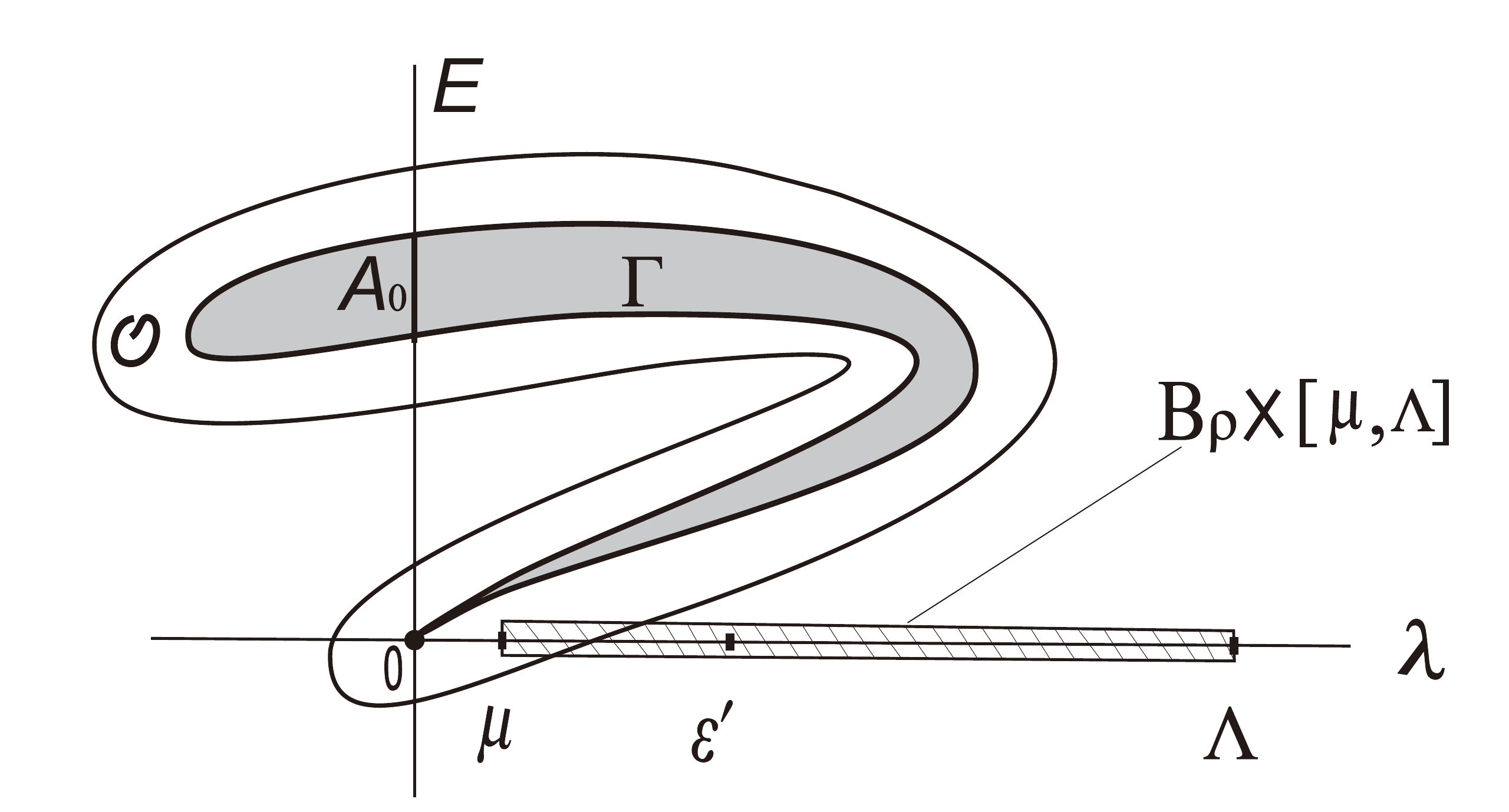}\includegraphics[width=6.6cm]{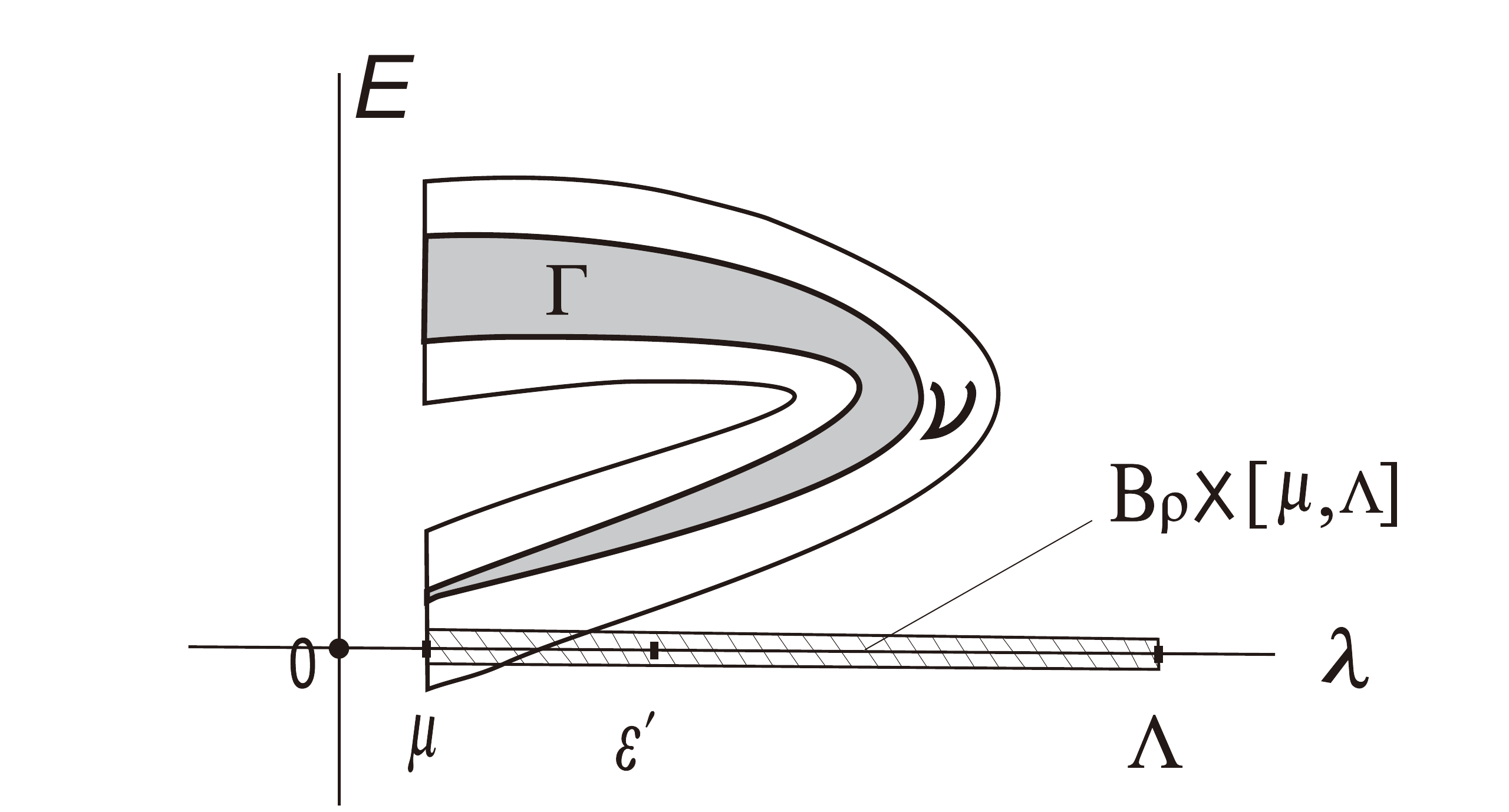}\\[2ex] Figure 6.7 \hspace{5cm} Figure 6.8 \end{center}

We first consider   $\lam\in[\mu,c]$. Because  $K_\lam$ contains all nonzero equilibria of $\Phi_\lam$ in $N$, in this case
one finds (by the choice of $\rho$) that there are no equilibrium points of $\Phi_\lam$ in $\mB_\rho$ other than $0$. Thus  by the  definition of $\cW$  we see  that  $e\in \cW[\lam]$.

 Now assume  $\lam>c$. We show that  $\cO[\lam]\cap \mB_\rho=\emp$, hence
$\cV[\lam]=\cW[\lam]=\cO[\lam]$, and the conclusion immediately follows. Let $u\in \cO[\lam]$. As $\cO\subset \mB_\sX(\Gamma,r)$,  there is a $(v,\lam')\in\Gamma$ such that
$$
||u-v||_\a+|\lam-\lam'|\leq r<\frac{1}{2}\min\{\kappa,c\},
$$
 which implies  $\lam'\geq c/2$ (otherwise $|\lam-\lam'|>c/2>r$). Hence by  \eqref{e6.8} we have $||v||_a\geq 4\kappa$. Therefore
$$
||u||_\a\geq ||v||_\a-||u-v||_\a\geq 4\kappa-r>3\kappa.
$$
We thus infer from the choice of $\rho$ that $u\not\in\mB_\rho$. This verifies that $\cO[\lam]\cap \mB_\rho=\emp$ and completes the proof of \eqref{5.16}.

\vs
It is trivial to check that  $\cV$ is a neighborhood of $\sC=\sC(\cV)=\sC(\cW)$ in $\sY=E\X[\mu,\Lam]$.
%Suppose the contrary. Then $\sC\cap \pa_{\sY}\cV\ne\emp$, where $\pa_\sY\cV$ denotes the boundary of $\cV$ relative to $\sY$. Observing  that $$\ba{ll}\pa_{\sY}\cV=\pa_{\sY}(\cO\cap \sY)\subset\pa\cO\cap \sY (??????),\ea$$ we have $$\ba{ll}\sC\cap \pa \cO=\sC\cap \(\pa\cO\cap \sY\)\supset \sC\cap \pa_{\sY}\cV \ne\emp.\ea$$ This contradicts  (\ref{5.10}). Note that $\sC$ is actually independent of the choice of $\rho$.
Pick a number  $0<\rho<\min(\kappa',\kappa)$ sufficiently small so that $\cW$ is a neighborhood of $\sC$ in $\sY$.
  Theorem \ref{t3.5} then asserts that
  \be\label{sp17} \mb{Ind}\,(\Phi_\lam,\sC[\lam])\equiv   \mb{Ind}(\Phi_{\Lam},\sC[\Lam])=  \mb{Ind}(\Phi_{\Lam},\emp)=0,\Hs \lam\in[\mu, \Lam].\ee
% \be\label{e:4.16b} \mb{Ind}\,(\Phi_\lam,\sC[\lam])\equiv   \mb{Ind}(\Phi_{\Lam},\sC_{\Lam})=  \mb{Ind}(\Phi_{\Lam},\emp)=0,\Hs \lam\in[\mu, \Lam].\ee
On the other hand, if $\mu\leq\lam\leq2\mu<c$ then
$$
\cW[\lam]=\cO[\lam]\sm\mB_\rho=(F_\lam\sm \mB_\rho)\cup G_\lam\cup H_\lam:=F_\lam'\cup G_\lam\cup H_\lam.
$$
Setting $$
\ell_0[\lam]=\sC[\lam]\cap F_\lam',\hs \ell_1[\lam]= \sC[\lam]\cap G_\lam,\hs \ell_2[\lam]=\sC[\lam]\cap H_\lam,$$
one obtains that
$$\mb{Ind}\,(\Phi_\lam,\sC[\lam])= \mb{Ind}(\Phi_{\lam},\ell_0[\lam])+ \mb{Ind}(\Phi_{\lam},\ell_1[\lam])+\mb{Ind}(\Phi_{\lam},\ell_2[\lam])$$
for $\lam\in[\mu, 2\mu]$. It follows  by \ef{sp17} that
\be\label{sp18}
 \mb{Ind}(\Phi_{\lam},\ell_1[\lam])+\mb{Ind}(\Phi_{\lam},\ell_2[\lam])=-\mb{Ind}(\Phi_{\lam},\ell_0[\lam]),\Hs \lam\in[\mu, 2\mu].
\ee

We infer from \eqref{e:4.15} and the choice of $\rho$ that  $K_\lam\subset F_\lam\sm \mB_\rho=F_\lam'$ for $0<\lam\leq2\mu$.  Recalling  that $K_\lam$ contains all nontrivial equilibria  in $N$ (and hence in $F_\lam'$) and noticing that $\ell_0[\lam]=\sE_{\Phi_\lam}(F_\lam')$,  we deduce that
$$\mb{Ind}(\Phi_{\lam},\ell_0[\lam])=\mb{Ind}(\Phi_{\lam},K_{\lam}),\Hs \lam\in [\mu,2\mu].$$
Thus by \ef{sp18} we have
\be\label{sp19}
\mb{Ind}(\Phi_{\lam},\ell_1[\lam])+\mb{Ind}(\Phi_{\lam},\ell_2[\lam])=-\mb{Ind}(\Phi_{\lam},K_{\lam}),\Hs \lam\in[\mu, 2\mu].
\ee
\vs
{\em Step 5.} Finally  we show that the left-hand side of \ef{sp19} equals $0$, thus obtain a contradiction.
For this purpose, consider the domain $\cG$  in $\sX_1=E\X(-\8,a]$  defined in \ef{spfn}. Let $\~\cG=\cG\cap \sW$, where $\sW=E\X[-\Lam,c]$, and $\Lam$ is the number in \ef{eol}. We claim that the boundary $\pa_{\sW}\~\cG$ of $\~\cG$ in $\sW$ contains no equilibrium points. To see this, by \ef{5.10} %Indeed, by \ef{sp9g} it is clear  that $(0,\lam)\not\in \cG\supset\~\cG$. Now we
it suffices to check  that $\pa_{\sW}\~\cG\subset \pa\cO$.

Indeed, since $\sW\subset \sX_1$, we infer from Lemma \ref{bl} that
$$
\pa_\sW(\~\cG)=\pa_\sW(\cG\cap\sW)\subset \pa_{\sX_1}(\cG).
$$
On the other hand, by \ef{cgh} we find that $\cG\subset \mB_{\sX_1}(\cR_a,\de_2)=\mb{int}_{\sX_1}Q,$ where $Q=\ol\mB_{\sX_1}(\cR_a,\de_2)$.
Therefore by Lemma
 \ref{bl} we deduce that 
 $$
 \pa_{\sX_1}(\cG)=\pa_{Q}(\cG)=\pa_{Q}(\cO\cap Q)\subset \pa \cO.
 $$
 Hence $\pa_{\sW}\~\cG\subset \pa\cO$.

Let $\wp=\sC(\~\cG)$. Then $\wp$ is a compact invariant set of the skew-product flow $\Pi$ of the family $\Phi_\lam$ ($\lam\in [-\Lam,c]$) on $\sW$. We infer from the claim proved above that $\wp\cap\pa_\sW\~\cG=\emp $. Hence $\wp$ is a static isolating set of $\Pi$. Thanks to Theorem  \ref{t3.5}, 
$$\mb{Ind}(\Phi_{\lam},\wp[\lam])\equiv \mb{Ind}(\Phi_{-\Lam},\wp[-\Lam])=\mb{Ind}(\Phi_{-\Lam},\emp)=0,\Hs \lam\in[-\Lam, c].$$
 Therefore 
$$
\mb{Ind}(\Phi_{\lam},\ell_1[\lam])=\mb{Ind}(\Phi_{\lam},\wp[\lam])= 0,\Hs \lam\in[\mu,2\mu].
$$

A parallel  argument as above applies to show that
$
\mb{Ind}(\Phi_{\lam},\ell_2[\lam])=0$ for $\lam\in[\mu,2\mu].$
Finally combining the above results with \ef{sp19} we conclude $\mb{Ind}(\Phi_{\lam},K_{\lam})=0$. This and \ef{sp7} contradict Theorem \ref{t5.1}. 

The  proof of the theorem is finished. \,$\Box$
\Vs
\br\label{gbr1}
If the third case (3) in Theorem \ref{gbt} occurs then  both $\Gamma^\pm$ are nontrivial; furthermore, we have  $d(0,\Gamma[0]\sm\{0\})>0$. Therefore as depicted in Fig. \ref{},  it is easy to see that there is a two-sided neighborhood $I_0$ of $\lam_0$ such that for each each $\lam\in I_0$ with $\lam\ne \lam_0$, the system $\Phi_\lam$ has at least two distinct nontrivial equilibria. Consequently we have a weaker version of Theorem \ref{gbt}:
\er
\bt\label{wgbt} Assume  the hypotheses in Theorem \ref{gbt}. Then either there is a two-sided neighborhood $I_0$ of $\lam_0$ such that for each $\lam\in I_0\sm\{\lam_0\}$,  $\Phi_\lam$ has at least two distinct nontrivial equilibria, or one of the following two assertions holds:
\begin{enumerate}
\item[$(1)$]  $\Gamma$ is unbounded.
%\item[(1)] There exist closed isolating neighborhood $N$ of $0$ and $\ve>0$ with $\cN:=N\X[-\ve,\ve]\subset \Om$ such that $\Gamma_\cN:=\Gamma_\cNM_0\subset \Gamma$; moreover,
% $$ \Gamma_\cN\cap (N\X\{-\ve\})\ne\emp\ne \Gamma_\cN\cap (N\X\{\ve\}).$$
\item[$(2)$] There exists $\lam_1\ne\lam_0$  such that $(0,\lam_1)\in\Gamma$\,.
\end{enumerate} \et

\section{The case $\mathfrak{m}_2=2$}
We now pay some attention to a particular case, namely, the case $\mathfrak{m}_2=2$. An easy example will also be included to illustrate our theoretical results.

\subsection{A local and global bifurcation theorem}

In what follows, by a {\em $k$-dimensional topological sphere} it means the boundary $\pa D$ of any contractible open subset $D$ of a $(k+1)$-dimensional manifold $\cM$ without boundary.
Denote $\mathbb{S}^k$   any $k$-dimensional topological sphere.

The main results in this section are summarized in the following theorem.

\bt\label{t8.1}Assume  (H1)-(H4)   are fulfilled with $\mathfrak{m}_2=2$. Suppose  $S_0=\{0\}$ is an isolated invariant set of $\Phi_{\lam_0}$ and $\chi\(h(\Phi_{\lam_0}^{2},S_0)\)\ne 1.$

Then one of the following two assertions holds.
\begin{enumerate}
\item[(1)]  There is a one-sided neighborhood $I_1$ of $\lam_0$ such that for each $\lam\in I_1\sm\{0\}$, the system $\Phi_\lam$ has a compact invariant set $S_\lam=\mathbb{S}^1$ with $0\not\in S_\lam$, and $S_\lam$  consists of either a closed orbit, or some nontrivial  equilibrium points of $\Phi_\lam$ and connecting orbits between them.
    \item[(2)]   $\Phi_\lam$ undergoes a static bifurcation as stated in Theorem \ref{wgbt}.
%\item[(2)] There is a two-sided neighborhood $I_2$ such that  $\Phi_\lam$ has at least two distinct nontrivial equilibria for each $\lam\in I_2\sm\{0\}$.
%\item[(3)]  The global static bifurcation branch $\Gamma$ of $(0,\lam_0)$ is either unbounded, or meets another bifurcation point $(0,\lam_1)\ne(0,\lam_1).$
    \end{enumerate}

\et

To prove the theorem, we need a basic result on the planar  system \be\label{e8.1}
\dot{x}=f(x),\Hs x=(x_1,x_2)\in\R^2.
\ee
Denote $\pi$ the local semiflow of  \ef{e8.1}.

Assume $f(0)=0$, and suppose $S_0=\{0\}$ is an isolated invariant set of $\pi$.
\bl\label{l8.2}
Suppose $S_0$ is neither an attractor nor a repeller. Then \be\label{eh2}H_q(h(\pi,S_0))=0,\Hs q\ne 1.\ee
\el
{\bf Proof.} It is known that $H_q(h(\pi,S_0))=0$ for $q>2$. So one only needs to verify the validity of \ef{eh2} for  $q=0,2$.

We infer from \cite[Theorem 1.5]{CE} that $S_0$ has an isolating block $B$  with smooth boundary $\pa B$. Note that $\pa B$ consists of Jordan curves.  Therefore there is at least one Jordan curve $C\subset \pa B$ such that $0\in\mb{int}\,N$, where $N$ denotes the bounded closed domain with $\pa N=C$. It is easy to understand that $N$ is an isolating block of $S_0$; moreover,  $N$ is contractible.

Because $S_0$ is neither an attractor nor a repeller, one has
$
  N^e\ne\emptyset\ne  N^i.
$
(See Subsection 2.3 for the definition of $N^e$, $N^i$ and $N^\pm$.) 
Thus we see that  $  N^-$ is the union of at most countably infinitely many disjoint curve segments $C_i$ ($1\leq i\leq n_0\leq\8$). For each $i$, we fix a point $p_i\in C_i$. Then  $p_i$ is a strong deformation retract of $C_i$. It follows that $P=\{p_i:\,\,1\leq i\leq n_0\}$ is a strong deformation retract of $  N^-$.
One can  easily  check that
\be\label{e8.2}
H_q(  N^-)=H_q(P)=0,\Hs \mb{for }q>0.
\ee

Recall that we have the exact sequence
$$
H_2(  N^-)\stackrel{i_*}{\longrightarrow}H_2(N)\stackrel{j_*}{\longrightarrow}H_2(N,  N^-)\stackrel{\pa_*}{\longrightarrow}H_1(  N^-).
$$
By \ef{e8.2} one concludes that $j_*$ is an isomorphism. %Hence $H_2(N,  N^-)\cong H_2(N)$. 
Therefore 
$$
H_2(h(\pi,S_0))\cong H_2(N,  N^-)\cong H_2(N)=0.
$$

Since $N$ is contractible, $N/  N^-$ is a path-connected space. Hence we have
$$
H_0(h(\pi,S_0))= H_0\((N/  N^-,[  N^-])\)=0,
$$
which completes the proof of the lemma. $\Box$

\Vs
Now we turn to the proof of Theorem \ref{t8.1}.
\Vs
\noindent{\bf Proof of Theorem \ref{t8.1}.} If $S_0$ is an attractor/repeller of $\Phi_{\lam_0}^2$, then the system undergoes an attractor/repeller bifurcation, and the conclusions in assertion  (1) follow from  the attractor bifurcation theory in Ma and Wang \cite{MW1} (see also \cite[Theorem 4.2]{LW}).
So we assume   $S_0$ is neither an attractor nor a repeller of $\Phi_{\lam_0}^2$.

Let $\b_q=\mb{rank}\,\(H_q(h(\Phi_{\lam_0}^{2},S_0))\)$ be the $q$-th Betti number of $h(\Phi_{\lam_0}^{2},S_0)$. By Lemma \ref{l8.2} we have $\b_q=0$ for all $q\ne 1$. Hence
$$
\chi\(h(\Phi_{\lam_0}^{2},S_0)\)=\Sig_{q=0}^\8(-1)^q \b_q=-\b_1\leq 0\ne 1.
$$
By Theorem \ref{wgbt} one concludes that assertion (2)  holds. $\Box$

\subsection{An example}
Finally let us  give a simple example to illustrate our theoretical results.

Consider the periodic problem on  $J=[-\pi,\pi]$:
\be\label{PB1}
\left\{\ba{ll}
-u''=\lam u+a(x)u^2+h(x,u),\hs x\in (-\pi,\pi);\\[1ex]
u(-\pi)=u(\pi),\,\,\,\,u'(-\pi)=u'(\pi),\ea\right.
\ee
where $a\in C^2(J)$, and $h\in C^2(J\X\R)$. Moreover,
\be\label{e8.4}
h(x,s)=O(|s|^3)\hs\mb{as }s\ra 0
\ee
uniformly for $x\in J$.

Let $X=L^2(J)$. Define an operator $A$ on $X$ to be the differential operator $-\frac{d^2}{dx^2}$ associated with the periodic boundary condition in \eqref{PB1}. Then
$$
\sig(A)=\sig_p(A)=\{\lam_n\}_{n=0}^\8,\hs\mb{where }\,\lam_n=n^2.
$$
The first eigenvalue $\lam_0$ is simple with an eigenfunction $e^O(x)=1/\sqrt{2\pi}$, and all the others are of multiplicity $2$. For $n\geq 1$, $A$ has a pair of eigenfunctions
$$
e^n_1=\frac{1}{\sqrt{\pi}}\sin nx,\hs e^n_2=\frac{1}{\sqrt{\pi}}\cos nx
$$
pertaining to $\lam_n$. 
The system $\{e^0,e^1_1,e^1_2,\cdots,e^n_1,e^n_2,\cdots\}$ forms a normal orthogonal basis of $X$.

Fix a number $\a\in[0,1)$, and let $E=X^\a$. Denote $||\.||_\a$ and $||\.||$ the norms on $E$ and $X$, respectively. Define $g:E\ra X$ as
$$
g(u)(x)=a(x)u^2(x)+h(x,u(x))\,\,\,(x\in J),\Hs \A\,u\in E.
$$
Then \eqref{PB1} can be written in an abstract form
\be\label{PB2}
Au=\lam u+g(u),\Hs u\in E.
\ee
Now we turn to    the bifurcation  of   \eqref{PB2} at each eigenvalue  $\lam_k$, $k\geq 1$.

Consider the corresponding evolution equation
\be\label{EPB}
u_t+L_\lam u=g(u),\Hs u\in E,
\ee
where $L_\lam=A-\lam I$.
We fix $k\geq1$ and, for simplicity, rewrite $$e^k_1=e_1,\hs e^k_2=e_2.$$ Denote $X_1=\mb{span}\{e_1,e_2\}$, and let
$X_2=X_1^\perp$. Then $X=X_1\oplus X_2$. Set
$$
E_i=E\cap X_i,\Hs i=1,2.
$$
Then  $E_1=X_1$,  and the norms $||\.||_\a$ and $||\.||$ are equivalent  on $E_1$.  Theorem \ref{rtt0} asserts that there is a convex neighborhood $V_1$ of $0$ in $E_1$ as well as a continuously differentiable mapping $\zeta:V_1\ra E_2$ such that $\cM^c=\{u_1+\zeta(u_1):\,\,u_1\in V_1\}$ is a local center manifold of \eqref{EPB} at $\lam=\lam_k$, and the reduction  equation  on $\cM^c$ reads
\be\label{e8.7}
\frac{du_1}{dt}=g(u_1+\zeta(u_1)),\Hs u_1\in V_1.
\ee

The mapping $\zeta$ can be approximated by some simpler ones. Indeed, we infer from  \cite[Chap. II, Theorem 2.3]{Ryba}  that if $\phi:V_1\ra E_2$ is a $C^1$-mapping with Lipschitz derivative $\phi'$ such that $\phi(V_1)\subset D(A)$ and
\be\label{e8.8}
||\De(u_1)||\leq M||u_1||^\b,\Hs\A\,u_1\in V_1
\ee
for some constants  $M>0$ and $\b>1$, where
$$
\De(u_1)=\phi'(u_1)[Lu_1-P^1g(u_1+\phi(u_1))]-[L\phi(u_1)-P^2g(u_1+\phi(u_1))]
$$
with $L=L_{\lam_k}$, then
$$
||\phi(u_1)-\zeta(u_1)||_\a\leq \~M||u_1||^\b,\Hs \A\,u_1\in V_1.
$$

We observe that if $\phi'(u_1)=O(||u_1||)$ as $||u_1||\ra 0$, then since $LE_1=\{0\}$,
\be\label{e8.9}\ba{ll}
\phi'(u_1)[Lu_1-P^1g(u_1+\phi(u_1))]&=-\phi'(u_1)[P^1g(u_1+\phi(u_1))]\\[1ex]
&=O(||u_1||^3),\hs\mb{as $||u_1||\ra 0$.}\ea
\ee
 For every $u_1=c_1e_1+c_2e_2\in V_1$, we also have
$$
\ba{ll}
P^2g(u_1+\phi(u_1))&=g(u_1+\phi(u_1))-P^1g(u_1+\phi(u_1))\\[1ex]
&=a(u_1+\phi(u_1))^2-P^1[a(u_1+\phi(u_1))^2]+O(||u_1||^3)\\[1ex]
&=c_1^2[ae_1^2-\<a e_1^2,e_1\>e_1-\<ae_1^2,e_2\>e_2]\\[1ex]
&\,\,\,\,\,+2c_1c_2 [ae_1e_2-\<a e_1e_2,e_1\>e_1-\<ae_1e_2,e_2\>e_2]\\[1ex]
&\,\,\,\,\,+c_2^2 [ae_2^2-\<a e_2^2,e_1\>e_1-\<ae_2^2,e_2\>e_2]+O(||u_1||^3)\\[1ex]
&:=c_1^2w_1+2c_1c_2w_0+c_2^2w_2+O(||u_1||^3).
\ea
$$
Note that $w_i\in X_2$ ($i=0,1,2$). As $$L_2=L|_{X_2}:D(L_2)=D(L)\cap X_2\ra X_2$$ is a one-one mapping, there exist $v_i\in D(L_2)\subset E_2$ such that
$$
Lv_i=L_2v_i=w_i,\Hs i=0,1,2.
$$
Thus if we define $\phi(u_1)$ as
$$
\phi(u_1)=c_1^2v_1+2c_1c_2v_0+c_2^2v_2,\Hs \A\, u_1=c_1e_1+c_2e_2\in V_1,
$$
then
$$
\ba{ll}
&L\phi(u_1)-P^2g(u_1+\phi(u_1))\\[1ex]
=&c_1^2Lv_1+2c_1c_2Lv_0+c_2^2Lv_2-P^2g(u_1+\phi(u_1))\\[1ex]
=&c_1^2w_1+2c_1c_2w_0+c_2^2w_2-P^2g(u_1+\phi(u_1))=O(||u_1||^3).
\ea$$
Combining this with \eqref{e8.9}, one finds that \eqref{e8.8} is fulfilled with $\b=3$. Hence 
$$
||\phi(u_1)-\zeta(u_1)||_a=O(||u_1||^3).
$$

Let $u_1(t)=c_1(t)e_1+c_2(t)e_2$. Then the  reduction equation \eqref{e8.7} reads
\be\label{e8.10}
\frac{du_1}{dt}=\<au_1^2,e_1\>e_1+\<au_1^2,e_2\>e_2+O(||u_1||^3),
\ee
or equivalently
\be\label{e8.11}
\left\{\ba{ll}
c_1'(t)=B_1(c_1,c_2)+O(||u_1||^3);\\[1ex]
c_2'(t)=B_2(c_1,c_2)+O(||u_1||^3),\ea\right.
\ee
where
$$
B_i=\<ae_1^2,e_i\>c_1^2+2\<ae_1e_2,e_i\>c_1c_2+\<ae_2^2,e_i\>c_2^2, \Hs i=1,2.
$$
Now we can state the following bifurcation result.

\bt\label{t8.3}Suppose that  the bilinear form $B_i$ is positive  definite for  $i=1$ or $2$, i.e., there exists $\gamma>0$ such that
$$B_i(c_1,c_2)\geq \gamma (c_1^2+c_2^2),\Hs\A\,(c_1,c_2\in\R^2).
$$
 Then one of the following alternatives  occurs.
\begin{enumerate}
\item[(1)]   There is a two-sided neighborhood $I_2$ such that for each $\lam\in I_2\sm\{0\}$, the problem \eqref{PB1} has at least two distinct nontrivial solutions.
\item[(2)]  The global static  bifurcation branch $\Gamma$ of $(0,\lam_k)$ is either unbounded, or meets another bifurcation point $(0,\lam_m)\ne(0,\lam_k).$
    \end{enumerate}
\et
{\bf Proof.} By the assumption of the theorem it is trivial to check  that $u_1=0$ is an isolated equilibrium of \eqref{e8.10}. Consequently $u=0$ is an isolated equilibrium of $\Phi_{\lam_k}$, where $\Phi_\lam$ denotes  the local semifow  of \eqref{EPB}. Since the system \eqref{EPB} has a Lyapunov function $V(u)$ which is precisely the variational functional of the problem \eqref{PB1},  we easily deduce that $S_0=\{0\}$ is  an isolating invariant set of $\Phi_{\lam_k}$.

By positivity of $B_i$ we see that
$$
c_i'(t)=B_i(c_1,c_2)+O(||u_1||^3)>0
$$
unless $c_1(t)=0=c_2(t)$,  which implies that $S_0$ is neither an attractor nor a repeller of \eqref{e8.10}. Now  the conclusion of the theorem immediately follows from Theorem \ref{t8.1}. $\Box$

\br A simple example in which the bilinear form $B_1$ is positive definite is the one where the function  $a(x)=\sin kx$ or $\cos kx$ for some $k\geq 0$. Indeed, if, say, $a(x)=\sin kx$, then
$$
\<ae_1^2,e_1\>=\pi^{-3/2}\int_{-\pi}^\pi \sin^4kx\mb{dx}>0, \hs \<ae_2^2,e_1\>=\pi^{-3/2}\int_{-\pi}^\pi \sin^2kx\cos^2kx\mb{dx}>0,
$$
and
$$\<ae_1e_2,e_1\>=\pi^{-3/2}\int_{-\pi}^\pi \sin^3kx\cos kx\mb{dx}=0,
$$
from which it is obvious that $B_1$ is positive definite.
\er

%%%%%%%%%%%%%%%%%%%%%%%%%%%%%%%%%%%%%%%%%%%%%%%%%%%%%%%%%%%%%%%%%%%%%%%%%%%
%%%%%%%%%%%%%%%%%%%%%%%%%%%%%%%%%%%%%%%%%%%%%%%%%%%%%%%%%%%%%%%%%%%%%%%%%%%%
\newpage
\centerline{\large\bf{ Appendix A:} Isomorphisms Induced by Projections}

\Vs\Vs
Let $X_\lam^i,$ $X_\lam^{ij}$, $P_\lam^i$ and  $P_\lam^{ij}$ be the same as in Section 4.1.
Since $P^3_\lam=I-(P^1_\lam+P^2_\lam)$, the continuity of $P^1_\lam$ and $P^2_\lam$ implies that $P_\lam^3$ is continuous in $\lam$ as well.

By (H3) we can   assume $J_0$ is chosen sufficiently small so that
$$
||P_\lam^i-P_{\lam_0}^i||\leq c<1/2,\Hs \A\,\lam\in J_0,\,\,i=1,2,3.\eqno(A1)
$$
Then
$$
||P_\lam^{ij}-P^{ij}_{\lam_0}||\leq 2c<1.
$$

As before,  we drop the subscript ``$\lam_0$'' and
rewrite $$X^i=X_{\lam_0}^i,\hs X^{ij}= X_{\lam_0}^{ij},\hs  P^i=P_{\lam_0}^i,\hs P^{ij}= P_{\lam_0}^{ij}.$$

\noindent{\bf Proposition A1.}\, {\em  For each $i=1,2,3$, the restriction  $P^i|_{X^i_\lam}$ of $P^i$ on $X_\lam^i$ is an isomorphism between $X^i_\lam$ and $X^i$. }
\Vs
\noindent{\bf Proof.} To prove Pro. A1,  let us first  verify  that $P^i|_{X^i_\lam}$ are  one-to-one mappings.

As $P_\lam^3=I-P_\lam^{12}$,  we deduce that
$$
||P_\lam^{3}-P^{3}||=||P_\lam^{12}-P^{12}||\leq 2c<1.\eqno(A2)
$$
In what follows  we argue by contradiction and suppose  $P^i|_{X^i_\lam}$ fails to be a  one-to-one mapping for some $i$. Then there would exist   $x_i\in X^i_\lam$ with $x_i\ne0$ such that $P^i x_i=0$. Further   by (A1) and (A2) we see that
$$
||x_i||=||P_\lam^ix_i||=||P_\lam^ix_i-P^ix_i||\leq 2c||x_i||<||x_i||,$$
a contradiction\,!

Now we show that $P^i|_{X^i_\lam}$ are isomorphisms. Since $P^i|_{X^i_\lam}$ are  one-to-one mappings,  by \ef{fie3}  one immediately concludes that $P^i|_{X^i_\lam}$ are isomorphisms for $i=1,2$. So we only need to consider the case where $i=3$.

Let $Q=P^3+P_\lam^{12}$. Then $$Q|_{X^3_\lam}=P^3|_{X^3_\lam}+P_\lam^{12}|_{X^3_\lam}=P^3|_{X^3_\lam}\,.$$
Because
$$Q=(I-P^{12})+P_\lam^{12}=I-(P^{12}-P_\lam^{12}),
$$
and  $||P^{12}-P_\lam^{12}||<1$, by the basic knowledge in linear functional analysis, we know that $Q:X\ra X$
is an isomorphism. To show that $P^3|_{X^3_\lam}$ is an isomorphism, there remains to check that $Q {X^3_\lam}=X^3$. For this purpose, it suffices to show that
$Q^{-1}X^3\subset X^3_\lam$.

We argue by contradiction and suppose the contrary. There would exist $u\not\in X^3_\lam$ such that $Qu\in X^3$. Let $u=x_\lam+x^3_\lam$, where $x_\lam\in X_\lam^{12}$, and $x_\lam^3\in X_\lam^3$. Then $x_\lam\ne 0$. We observe that
$$
Qu=(P^3+P_\lam^{12})u=P^3u+P_\lam^{12}(x_\lam+x^3_\lam)=x_\lam+P^3u\in X^3.
$$
Hence $x_\lam\in X^3$. Thereby we have $x_\lam\in X_\lam^{12}\cap X^3$.
 It follows that $$P_\lam^{12}x_\lam=x_\lam,\hs P^{12}x_\lam=0.$$
 Thus
 $$
 ||x_\lam||=||P_\lam^{12}x_\lam- P^{12}x_\lam||\leq c||x_\lam||<||x_\lam||.
 $$
 This leads to a contradiction and completes the proof of the proposition. $\Box$

%Now we show that $P^3|_{X^3_\lam}$  is an isomorphism. Let us  first verify that $P^3|_{X^3_\lam}$ is one-to-one. Suppose the contrary. Then there would exist   $w\in X^3_\lam$ with $w\ne0$ such that $P^3 w=0$. Hence
%$$
%w-P^{12}w=P^3 w=0,
%$$
%which implies $w=P^{12}w\in X^{12}$. Therefore $w\in X_\lam^3\cap X^{12}$. Thus by \ef{fie2},
%$$
%||w||=||P^{12}w||=||P^{12}_\lam w-P^{12}w||\leq 2c||w||<||w||,
%$$
%which leads to a contradiction! (The second equality in the above equation is due to that $w\in X_\lam^3$ and hence $P^{12}_\lam w=0$.)

\Vs
Now we define for each $\lam\in J_0$ a linear  operators $T_\lam$ on $X$ as follows:
$$
T_\lam u=\Sig_{1\leq j\leq 3}(P^j|_{X_\lam^j}P^j_\lam )\,u,\Hs u\in X.
$$
It is trivial to check that $T_\lam$ is an isomorphism with $T_{\lam_0}=I$. Clearly $T_\lam$ is continuous in $\lam$, and
%There is a  family of invertible bounded linear  operators $T=T_\lam$ on $X$  depending continuously on $\lam$, such that
$$
    T_\lam X^i_\lam=\Sig_{1\leq j\leq 3}(P^j|_{X_\lam^j}P^j_\lam)\, X^i_\lam=P^i|_{X_\lam^i} X^i_\lam= X^i,\Hs i=1,2,3.
$$
Thus we have
\Vs
\noindent{\bf Proposition A2.} {\em Under the assumptions (H1)-(H3), there exists  a  family of isomorphisms $T_\lam$ ($\lam\in J_0$) on $X$  depending continuously on $\lam$ with
$T_{\lam_0}=I$, such that $$
    T_\lam X^i_\lam=X^i_{\lam_0}:=X^i,\Hs i=1,2,3.\eqno(A3)
$$
}

\end{document}